\documentclass[onecolumn,11pt]{article}
\usepackage[top=1in, bottom=1in, left=1in, right=1in]{geometry}
\usepackage{amsmath, amsthm, amssymb}
\usepackage{graphicx}
\usepackage{epstopdf}
\usepackage{overpic}
\usepackage{cancel}
\usepackage{rotating}
\usepackage{url}
\usepackage{caption}
\usepackage{color}
\usepackage{rotating}
\usepackage{multirow}
\usepackage{wrapfig}
\usepackage{mathtools}

\usepackage{setspace}
\usepackage{palatino} % needs 10
\usepackage[bottom,flushmargin,hang,multiple]{footmisc}
\usepackage{lipsum}
\newcommand\blfootnote[1]{%
  \begingroup
  \renewcommand\thefootnote{}\footnote{#1}%
  \addtocounter{footnote}{-1}%
  \endgroup
}

\DeclareMathOperator*{\argmin}{arg\rm{}min}

\newcommand{\ba}{\mathbf{a}}

\newcommand{\bx}{\mathbf{x}}

\newcommand{\by}{\mathbf{y}}
\newcommand{\bz}{\mathbf{z}}

\newcommand{\bV}{\mathbf{V}}

\newcommand{\bX}{\mathbf{X}}

\newcommand{\bmu}{\boldsymbol{\mu}}
\newcommand{\bPsi}{\boldsymbol{\Psi}}

\newcommand{\bXi}{\boldsymbol{\Xi}}
\newcommand{\bxi}{\boldsymbol{\xi}}
\newcommand{\bSigma}{\boldsymbol{\Sigma}}
\newcommand{\bTheta}{\boldsymbol{\Theta}}

\usepackage{textcomp,xcolor}
\definecolor{MatlabCellColour}{RGB}{250,250,250}
\definecolor{MatPurp}{rgb}{.625,.1406,.9375}
\usepackage{listings}

\lstdefinestyle{customc}{
  belowcaptionskip=.25\baselineskip,
  breaklines=true,
  frame=L,
  xleftmargin=\parindent,
  language=Matlab,
  showstringspaces=false,
  basicstyle=\small\ttfamily,
  keywordstyle=\bfseries\color{white!30!black},
  identifierstyle=\color{blue},  
  commentstyle=\itshape\color{green!60!black},
  stringstyle=\color{MatPurp},
  backgroundcolor=\color{MatlabCellColour}
 }
 
\definecolor{blue}{rgb}{0,0,1}
\definecolor{darkgreen}{rgb}{0,0.5,0}
\definecolor{red}{rgb}{1,0,0}
\definecolor{teal}{rgb}{0,0.5,0.7}

\graphicspath{{Figures/}}

\title{\huge{ Discovering governing equations from data:} \\ \vspace{.1in} 
\Large{Sparse identification of nonlinear dynamical systems}}
\author{Steven L. Brunton$^{1*}$, Joshua L. Proctor$^2$, J. Nathan Kutz$^3$\\
\footnotesize{$^1$ Department of Mechanical Engineering, University of Washington, Seattle, WA 98195, United States}\\
\footnotesize{$^2$Institute for Disease Modeling, Bellevue, WA 98004, United States}\\ 
\footnotesize{$^3$ Department of Applied Mathematics, University of Washington, Seattle, WA 98195, United States}\\
}
\date{}
\begin{document}
\maketitle
\blfootnote{$^*$ Corresponding author. Tel.: +1 (609)-921-6415.\\ {\indent\emph{E-mail address:} sbrunton@uw.edu (S.L. Brunton).}}
%%%%%%%%%%%%
%%% ABSTRACT
%%%%%%%%%%%%
\vspace{-.2in}
\begin{abstract}
The ability to discover physical laws and governing equations from data is one of humankind's greatest intellectual achievements.  
A quantitative understanding of dynamic constraints and balances in nature has facilitated rapid development of knowledge and enabled advanced technological achievements, including aircraft, combustion engines, satellites, and electrical power.  
In this work, we combine sparsity-promoting techniques and machine learning with nonlinear dynamical systems to discover governing physical equations from measurement data.  
The only assumption about the structure of the model is that there are only a few important terms that govern the dynamics, so that the equations are sparse in the space of possible functions; this assumption holds for many physical systems. 
In particular, we use sparse regression to determine the fewest terms in the dynamic governing equations required to accurately represent the data. 
The resulting models are parsimonious, balancing model complexity with descriptive ability while avoiding overfitting.  
We demonstrate the algorithm on a wide range of problems, from simple canonical systems, including linear and nonlinear oscillators and the chaotic Lorenz system, to the fluid vortex shedding behind an obstacle.  The fluid example illustrates the ability of this method to discover the underlying dynamics of a system that took experts in the community nearly 30 years to resolve.  We also show that this method generalizes to parameterized, time-varying, or externally forced systems.\\

%\vspace{0.05in}
\end{abstract}

%%%%%%%%%%%%
%%% INTRODUCTION
%%%%%%%%%%%%
\section{Introduction}\label{sec:introduction}
Extracting physical laws from data is a central challenge in many diverse areas of science and engineering.  
There are many critical data-driven problems, such as understanding cognition from neural recordings, inferring patterns in climate, determining stability of financial markets, predicting and suppressing the spread of disease, and controlling turbulence for greener transportation and energy.  With abundant data and elusive laws, it is likely that data-driven discovery of {dynamics} will continue to play an increasingly important role in these efforts.

Advances in machine learning~\cite{Jordan2015science} and data science~\cite{Marx2013nature,Khoury2014science} have promised a renaissance in the analysis and understanding of complex data, extracting patterns in vast multimodal data that is beyond the ability of humans to grasp.  
However, despite the rapid development of tools to understand static data based on statistical relationships, there has been slow progress in distilling physical models of dynamic processes from big data.  
This has limited the ability of data science models to extrapolate the dynamics beyond the attractor where they were sampled and constructed.

An analogy may be drawn with the discoveries of Kepler and Newton.  
Kepler, equipped with the most extensive and accurate planetary data of the era, developed a \emph{data-driven} model for the motion of the planets, resulting in his famous elliptic orbits.
However, this was an \emph{attractor} based view of the world, and it did not explain the fundamental dynamic relationships that give rise to planetary orbits, or provide a model for how these bodies react when perturbed.  
Newton, in contrast, discovered a dynamic relationship between momentum and energy that described the underlying processes responsible for these elliptic orbits.  
This {dynamic} model may be generalized to predict behavior in regimes where no data was collected.  Newton's model has proven remarkably robust for engineering design, making it possible to land a spacecraft on the moon, which would not have been possible using Kepler's model alone. 

A seminal breakthrough by Schmidt and Lipson~\cite{Bongard2007pnas,Schmidt2009science} has resulted in a new approach to determine the underlying structure of a nonlinear dynamical system from data.  
This method uses symbolic regression (i.e., genetic programming~\cite{koza1992genetic}) to find nonlinear differential equations, and it balances complexity of the model, measured in the number of terms, with model accuracy.
The resulting {model identification} realizes a long-sought goal of the physics and engineering communities to discover dynamical systems from data. 
However, symbolic regression is expensive, does not scale well to large systems of interest, and may be prone to overfitting unless care is taken to explicitly balance model complexity with predictive power.  
In~\cite{Schmidt2009science}, the Pareto front is used to find parsimonious models in a large family of candidate models.  

In this work, we re-envision the dynamical system discovery problem from an entirely new perspective of sparse regression~\cite{Tibshirani1996lasso,Hastie2009book,James2013book} and compressed sensing~\cite{Donoho2006ieeetit,Candes2006bieeetit,Candes2006cpam,Candes2006picm,Baraniuk2007ieeespm,Tropp:2007}.
In particular, we leverage the fact that most physical systems have only a few relevant terms that define the dynamics, making the governing equations \emph{sparse} in a high-dimensional nonlinear function space.  
Before the advent of compressive sampling, and related sparsity-promoting methods, determining the few non-zero terms in a nonlinear dynamical system would have involved a combinatorial brute-force search, meaning that the methods would not scale to larger problems with Moore's law.  
However, powerful new theory guarantees that the sparse solution may be determined with high-probability using convex methods that do scale favorably with problem size.  
The resulting nonlinear model identification inherently balances model complexity (i.e., sparsity of right hand side dynamics) with accuracy, and the underlying convex optimization algorithms ensure that the method will be applicable to large-scale problems.

The method described here shares some similarity to the recent dynamic mode decomposition (DMD), which is a linear dynamic regression~\cite{Rowley2009jfm,schmid:2010}.  
DMD is an example of an equation-free method~\cite{Kevrekidis2003cms}, since it only relies on measurement data, but not on knowledge of the governing equations.  
Recent advances in the extended DMD have developed rigorous connections between DMD built on nonlinear observable functions and the Koopman operator theory for nonlinear dynamical systems~\cite{Rowley2009jfm,Mezic2013arfm}.  
However, there is currently no theory for which nonlinear observable functions to use, so that assumptions must be made on the form of the dynamical system.  
In contrast, the method developed here results in a \emph{sparse, nonlinear} regression that automatically determines the relevant terms in the dynamical system.  
The trend to exploit sparsity in dynamical systems is recent but growing~\cite{Schaeffer2013pnas,Ozolicnvs2013pnas,mackey2014compressive,Brunton2014siads,Proctor2014epj,Bai2014aiaa}.  
In this work, promoting sparsity in the dynamics results in parsimonious natural laws.\\

%%%%%%%%%%%%
%%% BACKGROUND
%%%%%%%%%%%%
\section{Background}\label{sec:background}
There is a long and fruitful history of modeling dynamics from data, resulting in powerful techniques for system identification~\cite{ljung:book}.  
Many of these methods arose out of the need to understand complex flexible structures, such as the Hubble space telescope or the international space station.  
The resulting models have been widely applied in nearly every branch of engineering and applied mathematics, most notably for model-based feedback control.  
However, methods for system identification typically require assumptions on the form of the model, and most often result in linear dynamics, limiting their effectiveness to small amplitude transient perturbations around a fixed point of the dynamics~\cite{guckenheimer_holmes}.  

This work diverges from the seminal work on system identification, and instead builds on symbolic regression and sparse representation.  
In particular, symbolic regression is used to find nonlinear functions describing the relationships between variables and measured dynamics (i.e., time derivatives).  
Traditionally, model complexity is balanced with describing capability using parsimony arguments such as the Pareto front.  
Here, we use sparse representation to determine the relevant model terms in an efficient and scalable framework.

%%% SYMBOLIC REGRESSION
\subsection{Symbolic regression and machine learning}\label{sec:SymReg}
Symbolic regression involves the determination of a function that relates input--output data, and it may be viewed as a form of machine learning.  
Typically, the function is determined using genetic programming, which is an evolutionary algorithm that builds and tests candidate functions out of simple building blocks~\cite{koza1992genetic}.  
These functions are then modified according to a set of evolutionary rules and generations of functions are tested until a pre-determined accuracy is achieved.  

Recently, symbolic regression has been applied to data from \emph{dynamical} systems, and ordinary differential equations were discovered from measurement data~\cite{Schmidt2009science}.  
Because it is possible to overfit with symbolic regression and genetic programming, a parsimony constraint must be imposed, and in~\cite{Schmidt2009science}, they accept candidate equations that are at the Pareto front of complexity.

%%% SPARSE REPRESENTATION
\subsection{Sparse representation and compressive sensing}\label{sec:SpaRep}
In many regression problems, only a few terms in the regression are important, and a \emph{sparse feature selection} mechanism is required.  
For example, consider data measurements $\by\in\mathbb{R}^m$ that may be a linear combination of columns from a feature library $\bTheta\in\mathbb{R}^{m\times p}$; the linear combination of columns is given by entries of the vector $\bxi\in\mathbb{R}^p$ so that:
\begin{eqnarray}
\by=\bTheta \bxi.\label{Eq:Axb}
\end{eqnarray}
Performing a standard regression to solve for $\bxi$ will result in a solution with nonzero contributions in each element.  
However, if \emph{sparsity} of $\bxi$ is desired, so that most of the entries are zero, then it is possible to add an $L^1$ regularization term to the regression, resulting in the LASSO~\cite{Hastie2009book,James2013book,Tibshirani1996lasso}:
\begin{eqnarray}
\bxi = \argmin_{\bxi'}\|\bTheta\bxi'-\by\|_2 + \lambda\|\bxi'\|_1.\label{Eq:LASSO}
\end{eqnarray}
The parameter $\lambda$ weights the sparsity constraint.  This formulation is closely related to the compressive sensing framework, which allows for the sparse vector $\bxi$ to be determined from relatively few \emph{incoherent} random measurements~\cite{Donoho2006ieeetit,Candes2006bieeetit,Candes2006cpam,Candes2006picm,Baraniuk2007ieeespm,Tropp:2007}.  
The sparse solution $\bxi$ to Eq.~\ref{Eq:Axb} may also be used for sparse classification schemes, such as the sparse representation for classification (SRC)~\cite{Wright2009ieeetpami}.  
Importantly, the compressive sensing and sparse representation architectures are convex and scale well to large problems, as opposed to brute-force combinatorial alternatives.

%%%%%%%%%%%%
%%% NONLINEAR SYSTEM ID
%%%%%%%%%%%%
\section{Sparse identification of nonlinear dynamics (SINDy)}\label{Sec:SparseDynamics}
In this work, we are concerned with identifying the governing equations that underly a physical system based on data that may be realistically collected in simulations or experiments.  
Generically, we seek to represent the system as a nonlinear dynamical system
\begin{eqnarray}
\dot\bx(t) = {\bf f}(\bx(t)).%,t;\bmu).
\label{Eq:CTsystem}
\end{eqnarray}
The vector $\bx(t)=\begin{bmatrix} x_1(t) & x_2(t) & \cdots & x_n(t)\end{bmatrix}^T\in\mathbb{R}^n$ represents the state of the system at time $t$, and the nonlinear function ${\bf f}(\bx(t))$ represents the dynamic constraints that define the equations of motion of the system.  
In the following sections, we will generalize Eq.~\eqref{Eq:CTsystem} to allow the dynamics $\bf{f}$ to vary in time, and also with respect to a set of bifurcation parameters $\bmu\in\mathbb{R}^q$.  
%The state $\bx$ may represent the positions and velocities of a physical system, the states of a spatial or spectral discretization of a partial differential equation, or the mode coefficients associated with a dimensionally reduced representation of a system.

The key observation in this paper is that for many systems of interest, the function $\bf{f}$ often consists of only a few terms, making it sparse in the space of possible functions.  
For example, the Lorenz system in Eq.~\eqref{Eq:Lorenz} has very few terms in the space of polynomial functions.  
Recent advances in compressive sensing and sparse regression make this viewpoint of sparsity favorable, since it is now possible to determine \emph{which} right hand side terms are non-zero without performing a computationally intractable brute-force search.  

To determine the form of the function $\bf{f}$ from data, we collect a time-history of the state $\bx(t)$ and its derivative $\dot\bx(t)$ sampled at a number of instances in time $t_1, t_2, \cdots, t_m$.  
These data are then arranged into two large matrices:
\begin{subequations}
\begin{eqnarray}
\bX &=& \begin{bmatrix} \bx^T(t_1)\\ \bx^T(t_2) \\ \vdots \\ \bx^T(t_m)\end{bmatrix} 
= \overset{\text{\normalsize state}}{\left.\overrightarrow{\begin{bmatrix}
x_1(t_1) & x_2(t_1) & \cdots & x_n(t_1)\\
x_1(t_2) & x_2(t_2) & \cdots & x_n(t_2)\\
\vdots & \vdots & \ddots & \vdots \\
x_1(t_m) & x_2(t_m) & \cdots & x_n(t_m)
\end{bmatrix}}\right\downarrow}\begin{rotate}{270}\hspace{-.125in}time~~\end{rotate}\label{Eq:DataMatrix}\\
\nonumber\\
\nonumber\\
\dot\bX & =&  \begin{bmatrix} \dot\bx^T(t_1)\\ \dot\bx^T(t_2) \\ \vdots \\ \dot\bx^T(t_m)\end{bmatrix} 
= \begin{bmatrix}
\dot x_1(t_1) & \dot x_2(t_1) & \cdots & \dot x_n(t_1)\\
\dot x_1(t_2) & \dot x_2(t_2) & \cdots & \dot x_n(t_2)\\
\vdots & \vdots & \ddots & \vdots \\
\dot x_1(t_m) & \dot x_2(t_m) & \cdots & \dot x_n(t_m)
\end{bmatrix}.\\\nonumber
\end{eqnarray}
\end{subequations}
Next, we construct an augmented library $\bTheta(\bX)$ consisting of candidate nonlinear functions of the columns of $\bX$.  
For example, $\bTheta(\bX)$ may consist of constant, polynomial and trigonometric terms:
\begin{eqnarray}
\bTheta(\bX) = 
\begin{bmatrix} 
~~\vline&\vline & \vline & \vline & & \vline & \vline & \vline& \vline &  ~~ \\
~~\mathbf{1}&\bX & \bX^{P_2} & \bX^{P_3} & \cdots & \sin(\bX) & \cos(\bX) & \sin(2\bX) & \cos(2\bX) & \cdots ~~\\
~~\vline &\vline & \vline & \vline & & \vline &\vline &\vline & \vline &  ~~
\end{bmatrix}.\label{Eq:NonlinearLibrary}
\end{eqnarray}
Here, higher polynomials are denoted as $\bX^{P_2}, \bX^{P_3},$ etc.  
For example, $\bX^{P_2}$ denotes the quadratic nonlinearities in the state variable $\bx$, given by:
\begin{eqnarray}
%\OneColEqu{
\bX^{P_2} = \begin{bmatrix}
x_1^2(t_1) & x_1(t_1)x_2(t_1) & \cdots& x_2^2(t_1) & x_2(t_1)x_3(t_1) & \cdots & x_n^2(t_1)\\
x_1^2(t_2) & x_1(t_2)x_2(t_2) & \cdots& x_2^2(t_2) & x_2(t_2)x_3(t_2) & \cdots & x_n^2(t_2)\\
\vdots &\vdots &\ddots&\vdots&  \vdots & \ddots & \vdots\\
x_1^2(t_m) & x_1(t_m)x_2(t_m) & \cdots & x_2^2(t_m) & x_2(t_m)x_3(t_m) & \cdots & x_n^2(t_m)
\end{bmatrix}.
%}
\end{eqnarray}

\begin{figure}
\begin{center}
\includegraphics[width=\textwidth]{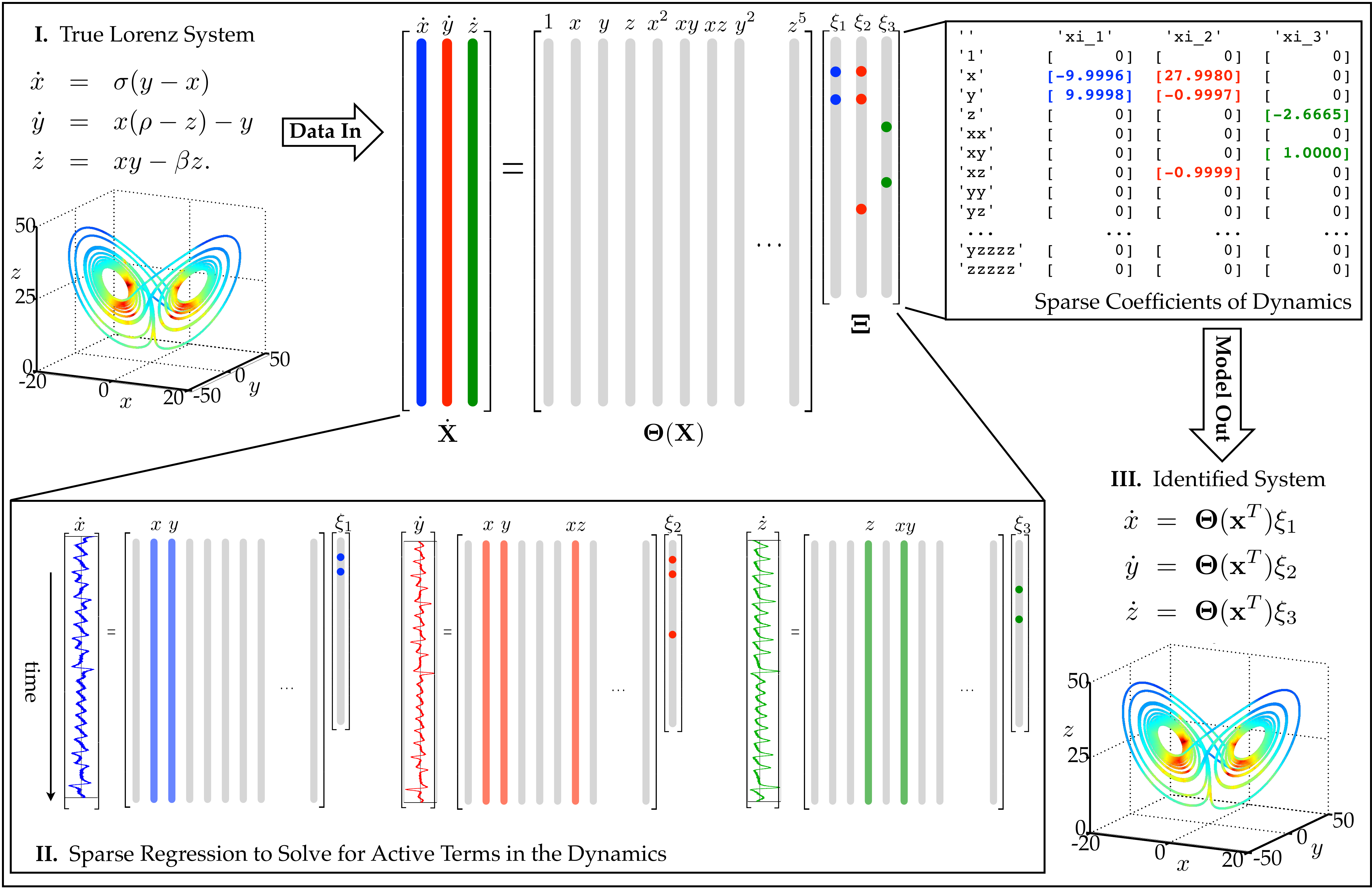}
\vspace{-.215in}
\caption{\small Schematic of our algorithm for sparse identification of nonlinear dynamics, demonstrated on the Lorenz equations.  Data is collected from measurements of the system, including {a time history of the states $\bX$ and derivatives $\dot{\bX}$}.  Next, a library of nonlinear functions of the states, $\bTheta(\bX)$, is constructed.  This nonlinear feature library is used to find the fewest terms needed to satisfy $\dot{\bX}=\bTheta(\bX)\bXi$.  The few entries in the vectors of $\bXi$, solved for by sparse regression, denote the relevant terms in the right-hand side of the dynamics.  Parameter values are $\sigma=10, \beta=8/3, \rho = 28$,  $(x_0,y_0,z_0)^T=(-8,7,27)^T$.  The trajectory on the Lorenz attractor is colored by the adaptive time-step required, with red requiring a smaller tilmestep.}\label{FIG00BIG}
\end{center}
\end{figure}

Each column of $\bTheta(\bX)$ represents a candidate function for the right hand side of Eq.~\eqref{Eq:CTsystem}.  
There is tremendous freedom of choice in constructing the entries in this matrix of nonlinearities.  
Since we believe that only a few of these nonlinearities are active in each row of $\bf{f}$, we may set up a sparse regression problem to determine the sparse vectors of coefficients $\bXi = \begin{bmatrix}\bxi_1 &\bxi_2 & \cdots & \bxi_n\end{bmatrix}$ that determine which nonlinearities are active, as illustrated in Fig.~\ref{FIG00BIG}.  
\begin{eqnarray}
\dot\bX = \bTheta(\bX) \bXi.\label{Eq:SparseRegression}
\end{eqnarray}

Each column $\bxi_k$ of $\bXi$ represents a sparse vector of coefficients determining which terms are active in the right hand side for one of the row equations $\dot{\bx}_k={\bf f}_k(\bx) $ in Eq.~\eqref{Eq:CTsystem}.  
Once $\bXi$ has been determined, a model of each row of the governing equations may be constructed as follows:
\begin{eqnarray}
\dot{\bx}_k = {\bf f}_k(\bx) = \bTheta(\bx^T)\bxi_k.\label{Eq:sparseRow}
\end{eqnarray}
Note that $\bTheta(\bx^T)$ is a vector of symbolic functions of elements of $\bx$, as opposed to $\bTheta(\bX)$, which is a data matrix.  This results in the overall model
\begin{eqnarray}
\dot{\bx} = {\bf f}(\bx) = \bXi^T (\bTheta(\bx^T))^T.\label{Eq:SparseFunction}
\end{eqnarray}
We may solve for $\bXi$ in Eq.~\eqref{Eq:SparseRegression} using sparse regression.

\subsection{Algorithm for sparse representation of dynamics with noise}\label{Sec:Algorithm}
There are a number of algorithms to determine sparse solutions $\bXi$ to the regression problem in Eq.~\eqref{Eq:SparseRegression}.  
Each column of Eq.~\eqref{Eq:SparseRegression} requires a distinct optimization problem to find the sparse vector of coefficients $\bxi_k$ for the $k^\text{th}$ row equation.  

For the examples in this paper, the matrix $\bTheta(\bX)$ has dimensions $m\times p$, where $p$ is the number of candidate nonlinear functions, and where $m\gg p$ since there are more time samples of data than there are candidate nonlinear functions.  
In most realistic cases, the data $\bX$ and $\dot\bX$ will be contaminated with noise so that Eq.~\eqref{Eq:SparseRegression} does not hold exactly.  
In the case that $\bX$ is relatively clean but the derivatives $\dot\bX$ are noisy, the equation becomes
\begin{eqnarray}
\dot\bX = \bTheta(\bX)\bXi + \eta \mathbf{Z},\label{Eq:NoisyRegression}
\end{eqnarray}
where $\mathbf{Z}$ is a matrix of independent identically distributed Gaussian entries with zero mean, and $\eta$ is the noise magnitude.  
Thus we seek a sparse solution to an overdetermined system with noise.  

%Common methods such as \texttt{cvx} and compressive sampling matching pursuit (CoSaMP)~\cite{Needell:2010} do not work well for these overdetermined systems.  
The LASSO~\cite{Hastie2009book,Tibshirani1996lasso} from statistics works well with this type of data, providing a sparse regression.  
However, it may be computationally expensive for very large data sets.  

An alternative is to implement the sequential thresholded least-squares algorithm in Code~\eqref{Code:sparseAlgorithm}.  
In this algorithm, we start with a least-squares solution for $\bXi$ and then threshold all coefficients that are smaller than some cutoff value $\lambda$.  
Once the indices of the remaining non-zero coefficients are identified, we obtain another least-squares solution for $\bXi$ onto the remaining indices.  
These new coefficients are again thresholded using $\lambda$, and the procedure is continued until the non-zero coefficients converge.  
This algorithm is computationally efficient, and it rapidly converges to a sparse solution in a small number of iterations.  
The algorithm also benefits from simplicity, with a single parameter $\lambda$ required to determine the degree of sparsity in $\bXi$.

Depending on the noise, it may still be necessary to filter the data $\bX$ and derivative $\dot\bX$ before solving for $\bXi$.  
In particular, if only the data $\bX$ is available, and $\dot\bX$ must be obtained by differentiation, then the resulting derivative matrix may have large noise magnitude.  
To counteract this, we use the total variation regularized derivative~\cite{Chartrand2011isrnam} to de-noise the derivative.  
An alternative would be to filter the data $\bX$ and $\dot\bX$, for example using the optimal hard threshold for singular values described in~\cite{Gavish2014arxiv}.

It is important to note that previous algorithms to identify dynamics from data have been quite sensitive to noise~\cite{Schmidt2009science}.  
The algorithm in Code~\eqref{Code:sparseAlgorithm} is remarkably robust to noise, and even when velocities must be approximated from noisy data, the algorithm works surprisingly well.

\begin{lstlisting}[caption={Sparse representation algorithm in Matlab.},label={Code:sparseAlgorithm}]
%% compute Sparse regression: sequential least squares
Xi = Theta\dXdt;  % initial guess: Least-squares

% lambda is our sparsification knob.
for k=1:10
    smallinds = (abs(Xi)<lambda);   % find small coefficients
    Xi(smallinds)=0;                % and threshold
    for ind = 1:n                   % n is state dimension
        biginds = ~smallinds(:,ind);
        % Regress dynamics onto remaining terms to find sparse Xi
        Xi(biginds,ind) = Theta(:,biginds)\dXdt(:,ind); 
    end
end
\end{lstlisting}

\subsection{Cross-validation to determine parsimonious sparse solution on Pareto front}
To determine the sparsification parameter $\lambda$ in the algorithm in Code~\eqref{Code:sparseAlgorithm}, it is helpful to use the concept of cross-validation from machine learning.  
It is always possible to hold back some test data apart from the training data to test the validity of models away from training values.  
In addition, it is important to consider the balance of model complexity (given by the number of nonzero coefficients in $\bXi$) with the model accuracy.  
There is an ``elbow" in the curve of accuracy vs. complexity parameterized by $\lambda$, the so-called Pareto front.  
This value of $\lambda$ represents a good tradeoff between complexity and accuracy, and it is similar to the approach taken in~\cite{Schmidt2009science}.

%\subsection{Systems with noise}

\subsection{Extensions and Connections}
There are a number of extensions to the basic theory above that generalize this approach to a broader set of problems.  First, the method is generalized to a discrete-time formulation, establishing a connection with the dynamic mode decomposition (DMD).  Next, high-dimensional systems obtained from discretized partial differential equations are considered, extending the method to incorporate dimensionality reduction techniques to handle big data.  Finally, the sparse regression framework is modified to include bifurcation parameters, time-dependence, and external forcing.  
 
%%% Discrete-time
\subsubsection{Discrete-time representation}
The aforementioned strategy may also be implemented on discrete-time dynamical systems:
\begin{eqnarray}
\bx_{k+1} &= &\mathbf{f}(\bx_k).
\label{Eq:DTsystem}
\end{eqnarray}
There are a number of reasons to implement Eq.~\eqref{Eq:DTsystem}.  First, many systems, such as the logistic map in Eq.~\eqref{Eq:logistic} are inherently discrete-time systems.  In addition, it may be possible to recover specific integration schemes used to advance Eq.~\eqref{Eq:CTsystem}.  The discrete-time formulation also foregoes the calculation of a derivative from noisy data.  The data collection will now involve two matrices $\bX_1^{m-1}$ and $\bX_2^m$:
\begin{eqnarray}
\bX_1^{m-1}=
\begin{bmatrix}\rule{.5in}{.75pt} & \bx_1^T &\rule{.5in}{.75pt} \\ 
\rule{.5in}{.75pt} & \bx_2^T &\rule{.5in}{.75pt} \\ 
& \vdots &  \\ 
\rule{.5in}{.75pt} & \bx_{m-1}^T &\rule{.5in}{.75pt}
\end{bmatrix}, \quad\quad
\bX_2^{m}=
\begin{bmatrix}\rule{.5in}{.75pt} & \bx_2^T &\rule{.5in}{.75pt} \\ 
\rule{.5in}{.75pt} & \bx_3^T &\rule{.5in}{.75pt} \\ 
& \vdots &  \\ 
\rule{.5in}{.75pt} & \bx_{m}^T &\rule{.5in}{.75pt}
\end{bmatrix}.
\end{eqnarray}
The continuous-time sparse regression problem in Eq.~\eqref{Eq:SparseRegression} now becomes:
\begin{eqnarray}
\bX_2^m = \bTheta(\bX_1^{m-1}) \bXi\label{Eq:DTSparseRegression}
\end{eqnarray}
and the function $\mathbf{f}$ is the same as in Eq.~\eqref{Eq:SparseFunction}.

In the discrete setting in Eq.~\eqref{Eq:DTsystem}, and for linear dynamics, there is a striking resemblance to dynamic mode decomposition.  In particular, if $\bTheta(\bx)=\bx$, so that the dynamical system is linear, then Eq.~\eqref{Eq:DTSparseRegression} becomes 
\begin{eqnarray}
\bX_2^m=\bX_1^{m-1}\bXi \quad\Longrightarrow \quad \left(\bX_2^m\right)^T=\bXi^T\left(\bX_1^{m-1}\right)^T.
\end{eqnarray} 
This is equivalent to the DMD, which seeks a dynamic regression onto linear dynamics $\bXi^T$.  In particular, $\bXi^T$ is $n\times n$ dimensional, which may be prohibitively large for a high-dimensional state $\bx$.  Thus, DMD identifies the dominant terms in the eigendecomposition of $\bXi^T$.

\subsubsection{High-dimensional systems, partial differential equations, and dimensionality reduction}
Often, the physical system of interest may be naturally represented by a partial differential equation (PDE) in a few spatial variables.  If data is collected from a numerical discretization or from experimental measurements on a spatial grid, then the state dimension $n$ may be prohibitively large.  For example, in fluid dynamics, even simple two-dimensional and three-dimensional flows may require tens of thousands up to billions of variables to represent the discretized system.  

The method described above is prohibitive for a large state dimension $n$, both because of the factorial growth of $\bTheta$ in $n$ and because each of the $n$ row equations in Eq.~\eqref{Eq:sparseRow} requires a separate optimization.  Fortunately, many high-dimensional systems of interest evolve on a low-dimensional manifold or attractor that may be well-approximated using a dimensionally reduced low-rank basis $\bPsi$~\cite{HLBR_turb,Majda2007bpnas}.  For example, if data $\bX$ is collected for a high-dimensional system as in Eq.~\eqref{Eq:DataMatrix}, it is possible to obtain a low-rank approximation using the singular value decomposition (SVD):
\begin{eqnarray}
\bX^T = \bPsi\bSigma\bV^*.
\end{eqnarray}
In this case, the state $\bx$ may be well approximated in a truncated modal basis $\bPsi_r$, given by the first $r$ columns of $\bPsi$ from the SVD:
\begin{eqnarray}
\bx & \approx &\bPsi_r \ba,
\end{eqnarray} 
where $\ba$ is an $r$-dimensional vector of mode coefficients.  We assume that this is a good approximation for a relatively low rank $r$.   
Thus, instead of using the original high-dimensional state $\bx$, it is possible to obtain a sparse representation of the Galerkin projected dynamics ${\bf f}_P$ in terms of the coefficients $\ba$:
\begin{eqnarray}
\dot\ba &=&{\bf f}_P(\ba).
\end{eqnarray}
There are many choices for a low-rank basis, including proper orthogonal decomposition (POD)~\cite{Berkooz:1993,HLBR_turb}, based on the SVD.%, balanced proper orthogonal decomposition (BPOD)~\cite{Willcox2002aiaaj,Rowley2005ijbc}, and dynamic mode decomposition (DMD)~\cite{Schmid2008aps,Rowley2009jfm,schmid:2010}, to name a few.  

\subsubsection{External forcing, bifurcation parameters, and normal forms}
In practice, many real-world systems depend on parameters, and dramatic changes, or bifurcations, may occur when the parameter is varied~~\cite{guckenheimer_holmes,Majda2009pnas}.  The algorithm above is readily extended to encompass these important parameterized systems, allowing for the discovery of normal forms associated with a bifurcation parameter $\bmu$.  First, we append $\bmu$ to the dynamics:
\begin{subequations}
\begin{eqnarray}
\dot{\mathbf{x}} & = &  \mathbf{f}(\mathbf{x};\bmu)\\
\dot\bmu & = & \mathbf{0}.
\end{eqnarray}
\end{subequations}
It is then possible to identify the right hand side $\mathbf{f}(\mathbf{x};\bmu)$ as a sparse combination of functions of components in $\mathbf{x}$ as well as the bifurcation parameter $\bmu$.  This idea is illustrated on two examples, the one-dimensional logistic map and the two-dimensional Hopf normal form. 

Time-dependence, known external forcing or control $\mathbf{u}(t)$ may also be added to the dynamics:
\begin{subequations}
\begin{eqnarray}
\dot{\mathbf{x}} & = &  \mathbf{f}(\mathbf{x},\mathbf{u}(t),t)\\
\dot t & = & 1.
\end{eqnarray}
\end{subequations}
This generalization makes it possible to analyze systems that are externally forced or controlled.  For example, the climate is both parameterized~\cite{Majda2009pnas} and has external forcing, including carbon dioxide and solar radiation.  The financial market presents another important example with forcing and active feedback control, in the form of regulations, taxes, and interest rates.

%%%%%%%%%%%%
%%% RESULTS
%%%%%%%%%%%%
\section{Results}
We demonstrate the methods  described in Sec.~\ref{Sec:SparseDynamics} on a number of canonical systems, ranging from simple linear and nonlinear damped oscillators, to noisy measurements of the fully chaotic Lorenz system, and to measurements of the unsteady fluid wake behind a circular cylinder, extending this method to nonlinear partial differential equations (PDEs) and high-dimensional data.
Finally, we show that bifurcation parameters may be included in the sparse models, recovering the correct normal forms from noisy measurements of the logistic map and the Hopf normal form.

\subsection{Example 1:  Simple illustrative systems}

\subsubsection{Example 1a: Two-dimensional damped oscillator (linear vs. nonlinear)}
\begin{figure}[where!]
\begin{center}
%\vspace{-.2in}
\begin{tabular}{cc}
Linear System & Cubic Nonlinearity\\
\begin{overpic}[width=.35\textwidth]{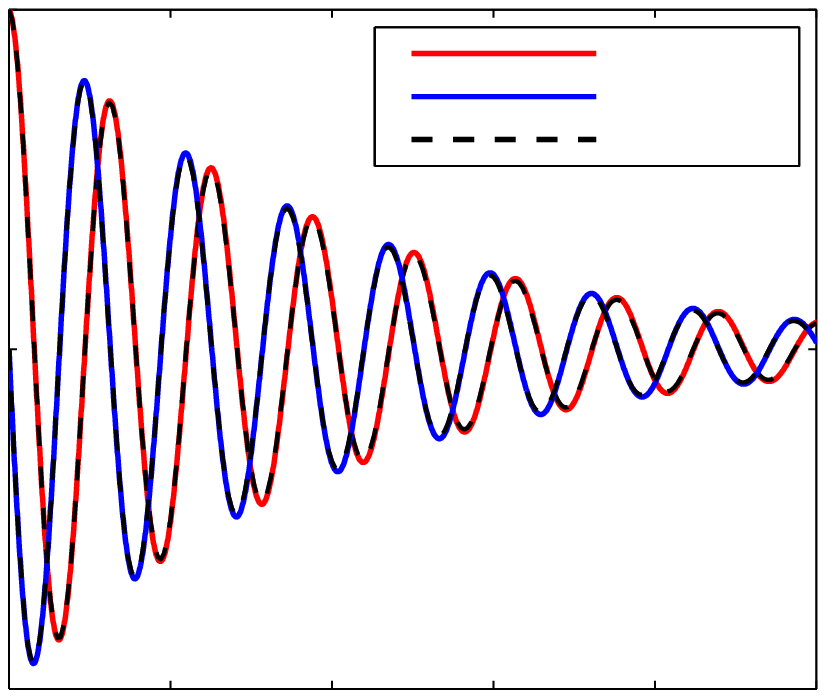}
{\small
\put(12,4){0}
\put(27.5,4){5}
\put(41.5,4){10}
\put(57,4){15}
\put(72.5,4){20}
\put(87.5,4){25}
%y-label
\put(8,8.5){-2}
\put(9.2,40){0}
\put(9.2,72){2}
% legend
\put(72,65){$x_2$}
\put(72,69){$x_1$}
\put(72,60.25){model}
}
\put(46,-2){Time}
\put(3.5,50){$x_k$}
\end{overpic}&
\begin{overpic}[width=0.35\textwidth]{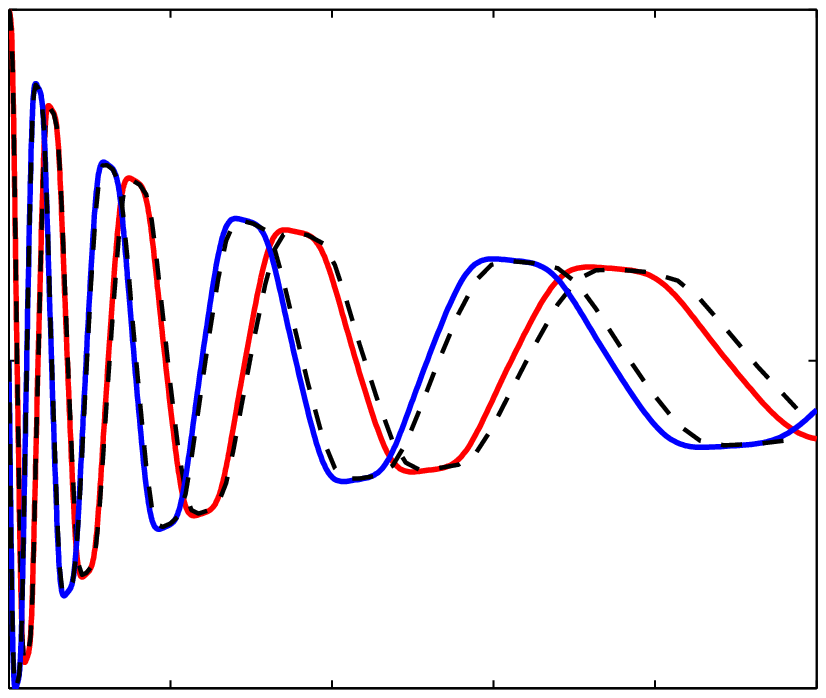}
{\small
\put(12,4){0}
\put(27.5,4){5}
\put(41.5,4){10}
\put(57,4){15}
\put(72.5,4){20}
\put(87.5,4){25}
%y-label
\put(8,8.5){-2}
\put(9.2,40){0}
\put(9.2,72){2}
}
\put(46,-2){Time}
\put(3.5,50){$x_k$}
\end{overpic}\\

%&\\
\begin{overpic}[width=0.35\textwidth]{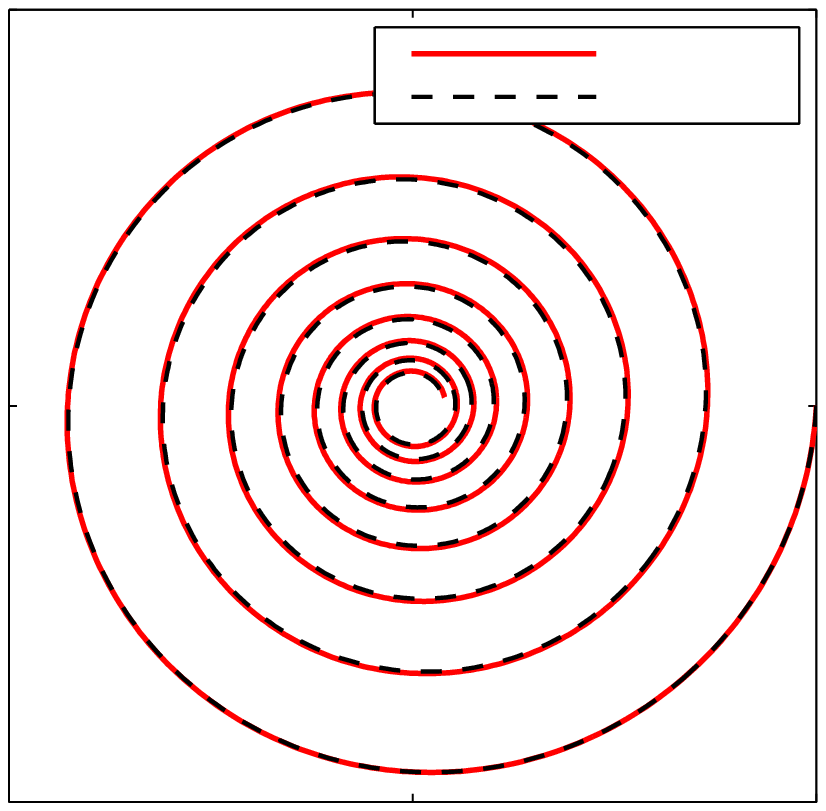}
{\small
\put(10,5.5){-2}
\put(50.5,5.5){0}
\put(89.5,5.5){2}
%y-label
\put(8,9.5){-2}
\put(9.2,48){0}
\put(9.2,85){2}
% legend
\put(72,81.25){$x_k$}
\put(72,77){model}
}
\put(55,0){$x_1$}
\put(3.5,55){$x_2$}
\end{overpic}&
\begin{overpic}[width=0.35\textwidth]{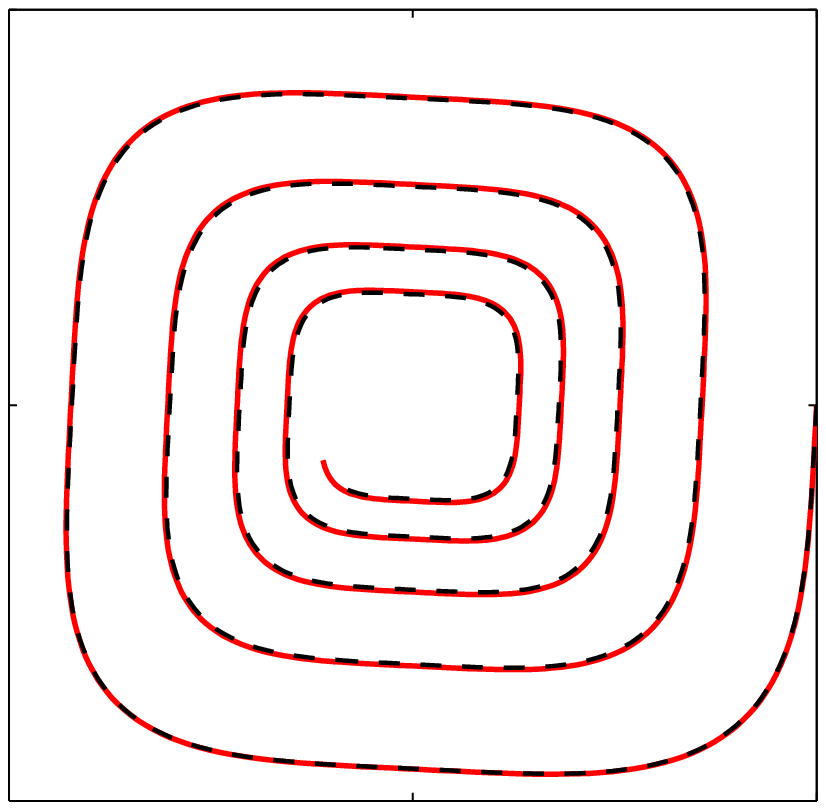}
{\small
\put(10,5.5){-2}
\put(50.5,5.5){0}
\put(89.5,5.5){2}
%y-label
\put(8,9.5){-2}
\put(9.2,48){0}
\put(9.2,85){2}
}
\put(55,0){$x_1$}
\put(3.5,55){$x_2$}
\end{overpic}
\end{tabular}
\vspace{-.1in}
\caption{\small Comparison of linear system (left) and system with cubic nonlinearity (right).   The sparse identified system correctly identifies the form of the dynamics and accurately reproduces the phase portraits.}\label{Fig:Ex1}
\end{center}
\vspace{-.1in}
\end{figure}

In this example, we consider the two-dimensional damped harmonic oscillator with either linear or cubic dynamics, as in Eq.~\eqref{Eq:CubicOsc}.  The dynamic data and the sparse identified model are shown in Fig.~\ref{Fig:Ex1}.   The correct form of the nonlinearity is obtained in each case; the augmented nonlinear library $\bTheta(\bx)$ includes polynomials in $\bx$ up to fifth order.  The sparse identified model and algorithm parameters are shown in the Appendix in Tables~\ref{Tab:Ex1_2dLin} and \ref{Tab:Ex1_2dCub}.
%\begin{subequations}
\begin{eqnarray}
\frac{d}{dt}\begin{bmatrix} x\\ y\end{bmatrix} =\begin{bmatrix} -0.1 & 2\\ -2 & -0.1\end{bmatrix}\begin{bmatrix} x \\ y \end{bmatrix}%%\\\label{Eq:LinOsc}
%\nonumber\\
\quad \quad \quad \quad
\frac{d}{dt}\begin{bmatrix} x\\ y\end{bmatrix}  =\begin{bmatrix} -0.1 & 2\\ -2 & -0.1\end{bmatrix}\begin{bmatrix} x^3 \\ y^3 \end{bmatrix}\label{Eq:CubicOsc}
\end{eqnarray}
%\end{subequations}

\subsubsection{Example 1b: Three-dimensional linear system}
\begin{figure}[where!]
\begin{center}
\begin{tabular}{cc}
\begin{overpic}[width=0.35\textwidth]{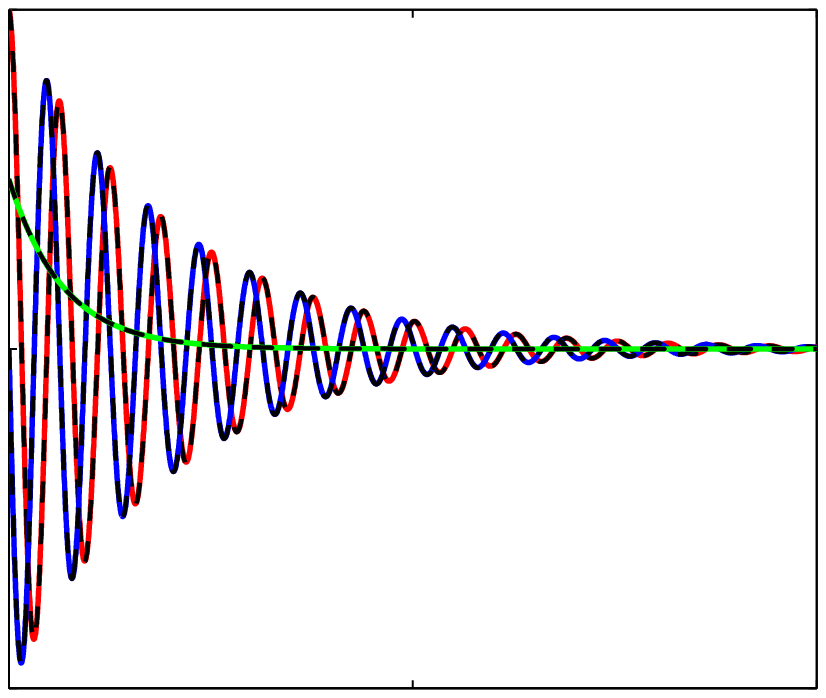}
{\small
\put(12,4){0}
\put(50,4){25}
\put(87.5,4){50}
%y-label
\put(8,8.5){-2}
\put(9.2,40){0}
\put(9.2,72){2}
}
\put(46,-2){Time}
\put(3.5,50){$x_k$}
\end{overpic}
&
\begin{overpic}[width=0.35\textwidth]{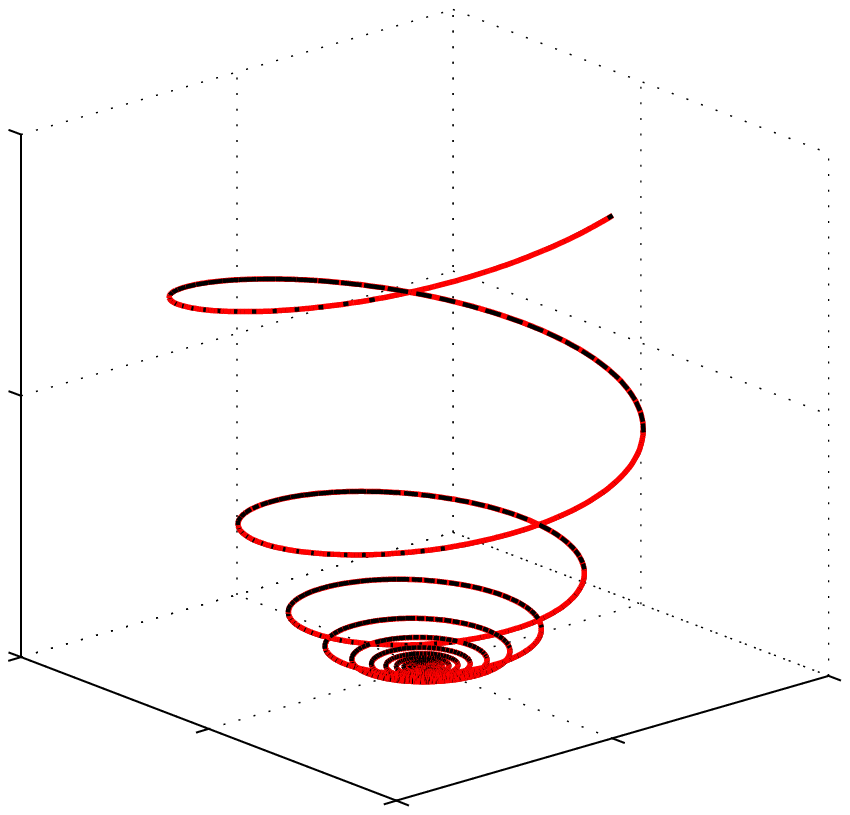}
{\small
%x-axis
\put(7,20.5){-2}
\put(26.5,13.5){0}
\put(44.5,7.5){2}
%y-axis
\put(49.5,7.){-2}
\put(71.5,13.){0}
\put(92,19.){2}
%z-axis
\put(8.5,24.5){0}
\put(5.75,49.){0.5}
\put(8.5,74.5){1}
}
\put(28,8){$x_1$}
\put(76,9){$x_2$}
\put(2,58){$x_3$}
\end{overpic}
\end{tabular}
\caption{\small Three-dimensional linear system (solid colored lines) is well-captured by sparse identified system (dashed black line).}\label{Fig:Ex1_3d}
\end{center}
\end{figure}

A linear system with three variables and the sparse approximation are shown in Fig.~\ref{Fig:Ex1_3d}.  In this case, the dynamics are given by
\begin{eqnarray}
\frac{d}{dt}\begin{bmatrix}x\\ y\\ z\end{bmatrix} & = & \begin{bmatrix}  -0.1 & -2  & 0 \\ 2 & -0.1 & 0 \\ 0 & 0 & -0.3\end{bmatrix}\begin{bmatrix} x\\ y \\ z\end{bmatrix}.
\end{eqnarray}
The sparse identification algorithm correctly identifies the system in the space of polynomials up to second or third order, and the sparse model is given in Table~\ref{Tab:Ex1_3d}.  Interestingly, including polynomial terms of higher order (e.g. orders 4 or 5) introduces a degeneracy in the sparse identification algorithm, because linear combinations of powers of $e^{\lambda t}$ may approximate other exponential rates.  This unexpected degeneracy motivates a hierarchical approach to identification, where subsequently higher order terms are included until the algorithm either converges or diverges.

%%%%%%%
\subsection{Example 2:  Lorenz system (Nonlinear ODE)}
Here, we consider the nonlinear Lorenz system~\cite{Lorenz1963jas} to explore the identification of chaotic dynamics evolving on an attractor, shown in Fig.~\ref{FIG00BIG}:
\begin{subequations}
\label{Eq:Lorenz}
\begin{eqnarray}
\dot{x} & = & \sigma (y - x)\\
\dot{y} & = & x(\rho -z) - y\\
\dot{z} & = & x y - \beta z.
\label{Eq:Lorenz}
\end{eqnarray}
\end{subequations}
Although these equations give rise to rich and chaotic dynamics that evolve on an attractor, there are only a few terms in the right-hand side of the equations.  Figure~\ref{FIG00BIG} shows a schematic of how data is collected for this example, and how sparse dynamics are identified in a space of possible right-hand side functions using convex $\ell_1$-minimzation.  

For this example, data is collected for the Lorenz system, and stacked into two large data matrices $\bX$ and $\dot\bX$, where each row of $\bX$ is a snapshot of the state $\bx$ in time, and each row of $\dot\bX$ is a snapshot of the time derivative of the state $\dot\bx$ in time.  Here, the right-hand side dynamics are identified in the space of polynomials $\bTheta(\bX)$ in $(x,y,z)$ up to fifth order:
\begin{eqnarray}
\bTheta(\bX)=\begin{bmatrix} 
\vline & \vline & \vline & \vline & \vline & \vline & \vline & \vline & \vline &  ~~~~~~~\vline \\
\bx(t) & \by(t) & \bz(t) & \bx(t)^2 & \bx(t)\by(t) & \bx(t)\bz(t) & \by(t)^2 & \by(t)\bz(t) & \bz(t)^2 & \cdots ~~~\bz(t)^5 \\
\vline & \vline & \vline & \vline & \vline & \vline & \vline & \vline & \vline & ~~~~~~~\vline\\
\end{bmatrix}.
\end{eqnarray}
Each column of $\bTheta(\bX)$ represents a candidate function for the right hand side of Eq.~\eqref{Eq:CTsystem}, and a sparse regression determines which terms are active in the dynamics, as in Fig.~\ref{FIG00BIG}, and Eq.~\eqref{Eq:SparseRegression}.

\begin{figure}[b!]
%\vspace{0.025in}
\begin{center}
\begin{tabular}{ccccccc}
\begin{overpic}[width=0.25\textwidth]{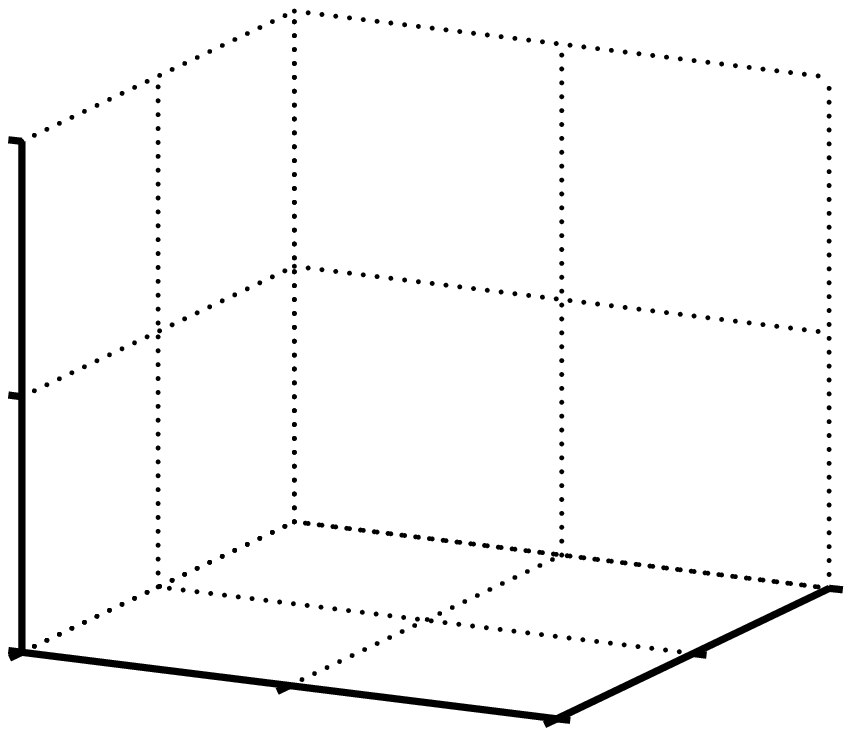}
\put(-12.5,-7.75){\begin{overpic}[width=.3\textwidth]{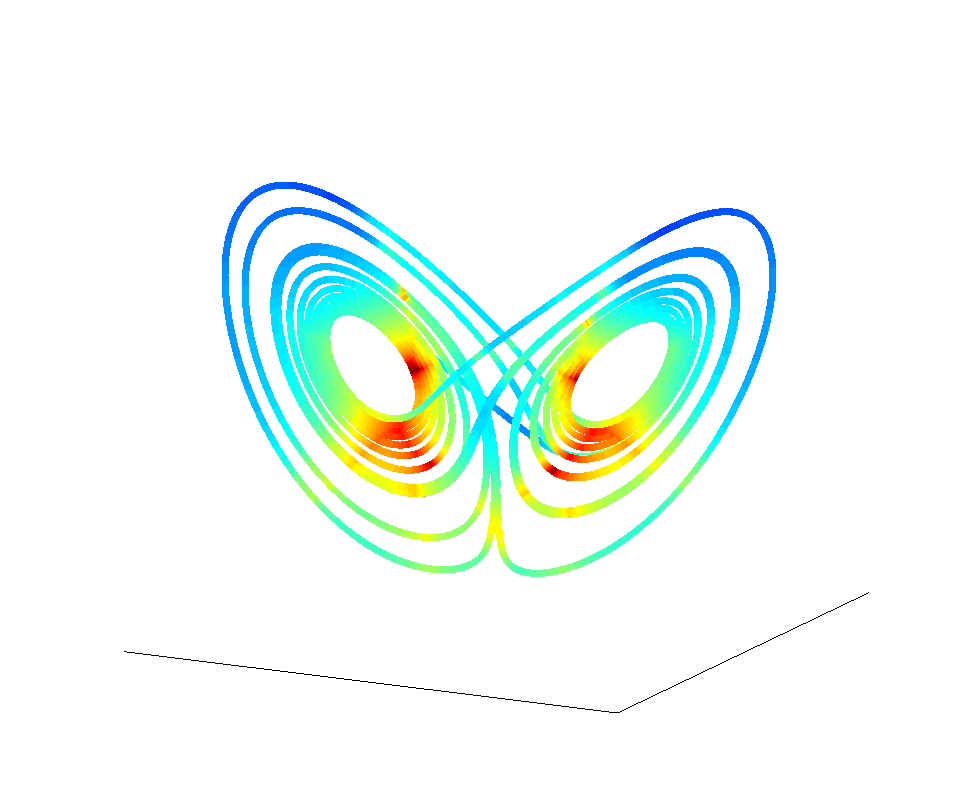}\end{overpic}}
\put(5,87){Full Simulation}
% label z-axis
{\small
\put(-3,9.5){0}
\put(-6,38){25}
\put(-6,67){50}}
\put(-8,48){$z$}
% label x-axis
{\small
\put(-4.5,3.5){-20}
\put(28,0){0}
\put(54.25,-3.){20}}
\put(36.5,-4.5){$x$}
% label y-axis
{\small
\put(66.5,-2.25){-50}
\put(83.5,6){0}
\put(97.5,12.5){50}}
\put(89,2){$y$}
\end{overpic} && &
\begin{overpic}[width=0.25\textwidth]{FIG02a_Lorenz_Axis.eps}
\put(-12.5,-7.75){\begin{overpic}[width=.3\textwidth]{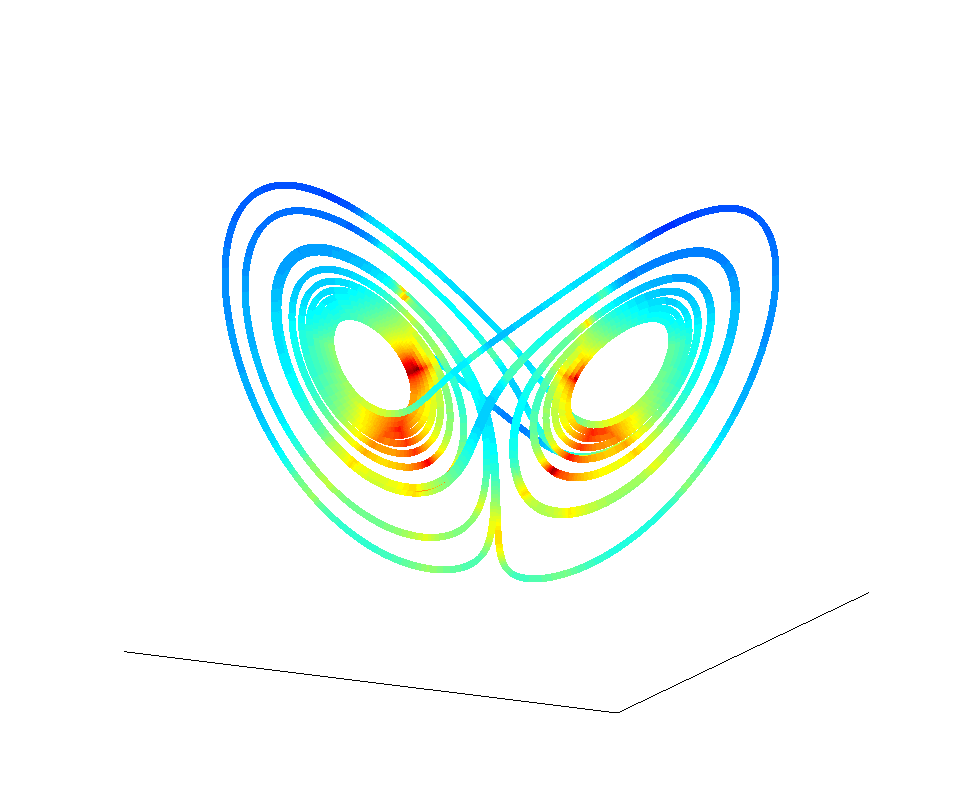}\end{overpic}}
\put(3,87){Identified System, $\eta=0.01$}
% label z-axis
{\small
\put(-3,9.5){0}
\put(-6,38){25}
\put(-6,67){50}}
\put(-8,48){$z$}
% label x-axis
{\small
\put(-4.5,3.5){-20}
\put(28,0){0}
\put(54.25,-3.){20}}
\put(36.5,-4.5){$x$}
% label y-axis
{\small
\put(66.5,-2.25){-50}
\put(83.5,6){0}
\put(97.5,12.5){50}}
\put(89,2){$y$}
\end{overpic}& &&
\begin{overpic}[width=0.25\textwidth]{FIG02a_Lorenz_Axis.eps}
\put(-12.5,-7.75){\begin{overpic}[width=.3\textwidth]{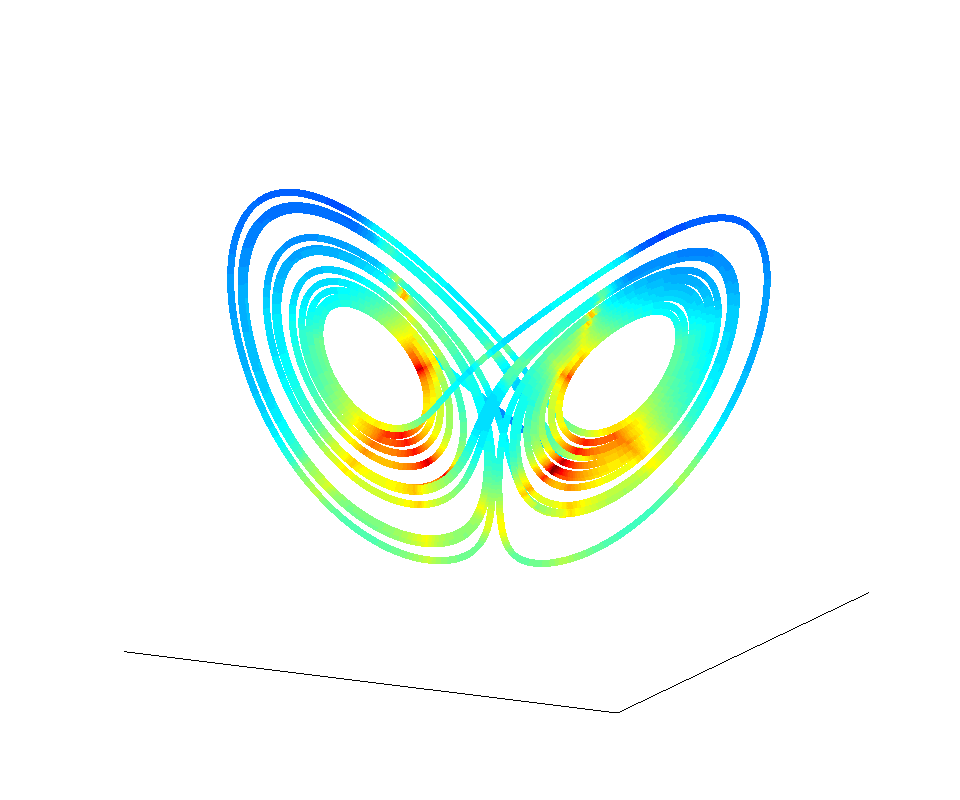}\end{overpic}}
\put(3,87){Identified System, $\eta=10$}
% label z-axis
{\small
\put(-3,9.5){0}
\put(-6,38){25}
\put(-6,67){50}}
\put(-8,48){$z$}
% label x-axis
{\small
\put(-4.5,3.5){-20}
\put(28,0){0}
\put(54.25,-3.){20}}
\put(36.5,-4.5){$x$}
% label y-axis
{\small
\put(66.5,-2.25){-50}
\put(83.5,6){0}
\put(97.5,12.5){50}}
\put(89,2){$y$}
\end{overpic} \\
&&\\
\begin{overpic}[width=0.25\textwidth]{FIG02a_Lorenz_Axis.eps}
\put(-12.5,-7.75){\begin{overpic}[width=.3\textwidth]{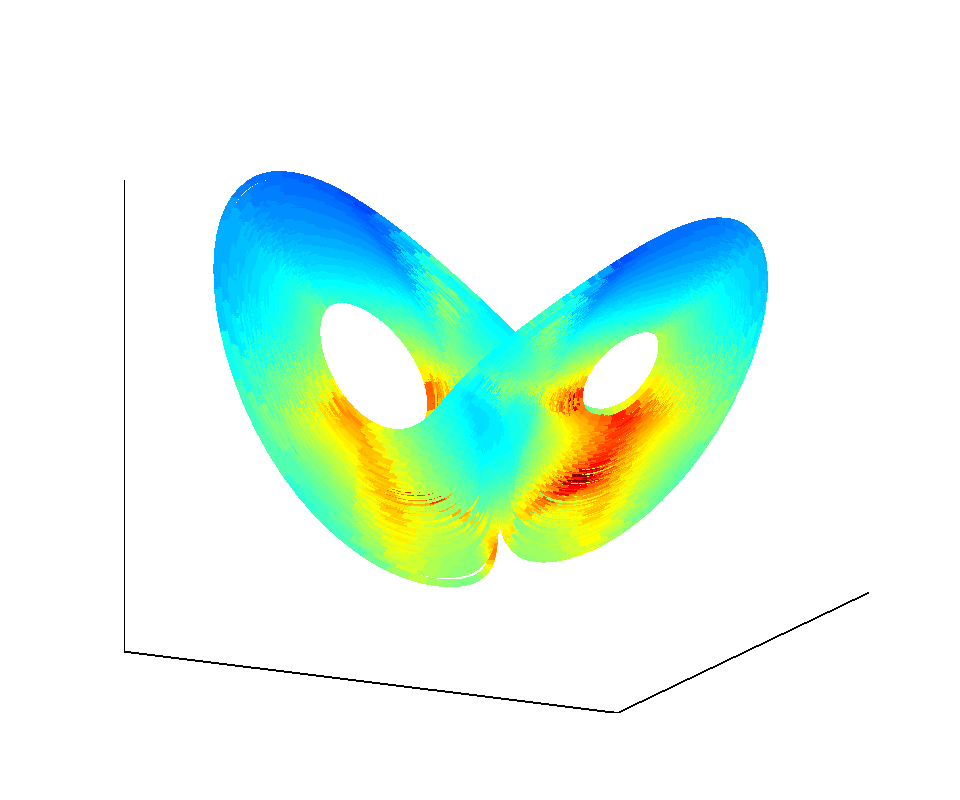}\end{overpic}}
%\put(5,87){Full Simulation}
% label z-axis
{\small
\put(-3,9.5){0}
\put(-6,38){25}
\put(-6,67){50}}
\put(-8,48){$z$}
% label x-axis
{\small
\put(-4.5,3.5){-20}
\put(28,0){0}
\put(54.25,-3.){20}}
\put(36.5,-4.5){$x$}
% label y-axis
{\small
\put(66.5,-2.25){-50}
\put(83.5,6){0}
\put(97.5,12.5){50}}
\put(89,2){$y$}
\end{overpic} && &
\begin{overpic}[width=0.25\textwidth]{FIG02a_Lorenz_Axis.eps}
\put(-12.5,-7.75){\begin{overpic}[width=.3\textwidth]{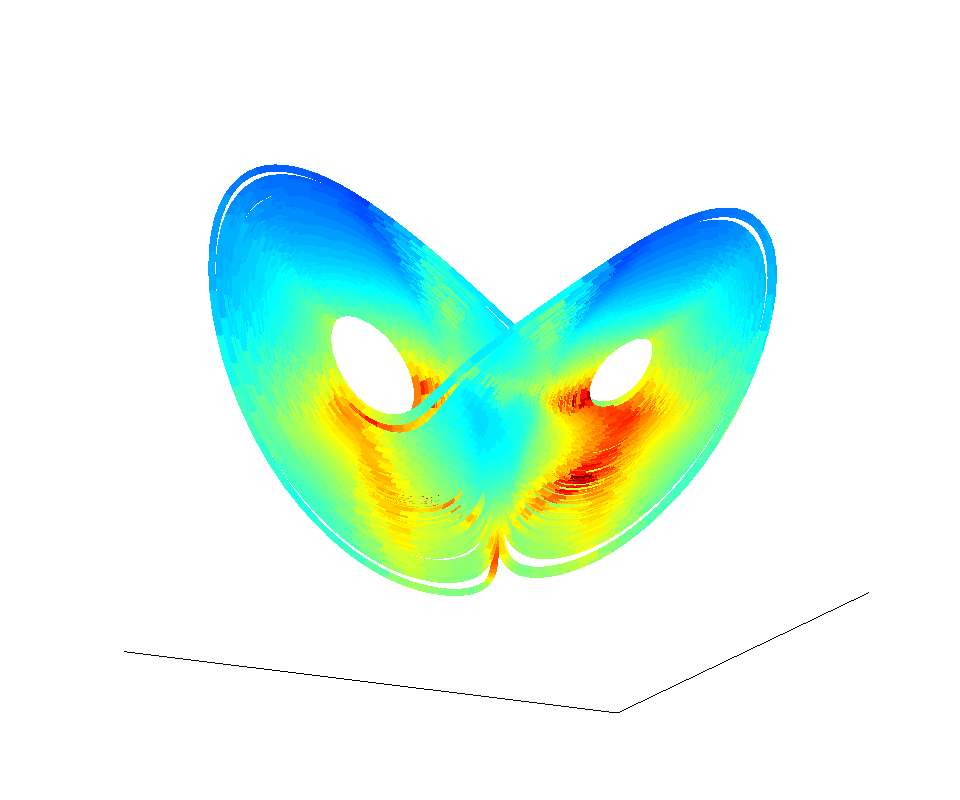}\end{overpic}}
%\put(5,87){Identified System, $\eta=0.01$}
% label z-axis
{\small
\put(-3,9.5){0}
\put(-6,38){25}
\put(-6,67){50}}
\put(-8,48){$z$}
% label x-axis
{\small
\put(-4.5,3.5){-20}
\put(28,0){0}
\put(54.25,-3.){20}}
\put(36.5,-4.5){$x$}
% label y-axis
{\small
\put(66.5,-2.25){-50}
\put(83.5,6){0}
\put(97.5,12.5){50}}
\put(89,2){$y$}
\end{overpic}& &&
\begin{overpic}[width=0.25\textwidth]{FIG02a_Lorenz_Axis.eps}
\put(-12.5,-7.75){\begin{overpic}[width=.3\textwidth]{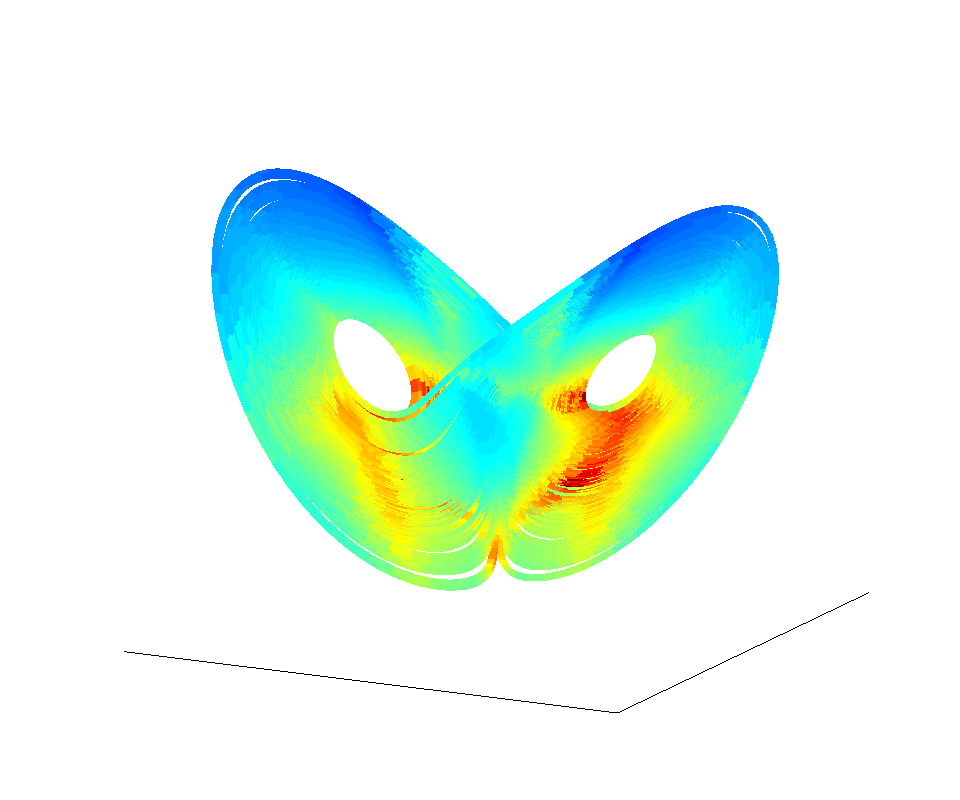}\end{overpic}}
%\put(5,87){Identified System, $\eta=10$}
% label z-axis
{\small
\put(-3,9.5){0}
\put(-6,38){25}
\put(-6,67){50}}
\put(-8,48){$z$}
% label x-axis
{\small
\put(-4.5,3.5){-20}
\put(28,0){0}
\put(54.25,-3.){20}}
\put(36.5,-4.5){$x$}
% label y-axis
{\small
\put(66.5,-2.25){-50}
\put(83.5,6){0}
\put(97.5,12.5){50}}
\put(89,2){$y$}
\end{overpic}
\end{tabular}
%\vspace{.1in}
\caption{\small Trajectories of the Lorenz system for short-time integration from $t=0$ to $t=20$ (top) and long-time integration from $t=0$ to $t=250$ (bottom).  The full dynamics (left) are compared with the sparse identified systems (middle, right) for various additive noise.  The trajectories are colored by $\Delta t$, the adaptive Runge-Kutta time step.  This color is a proxy for local sensitivity.}\label{Lorenz:Long}
\end{center}
\end{figure}

\begin{figure}[where!]
\begin{center}
\begin{tabular}{ccc}
\begin{overpic}[width=0.5\textwidth]{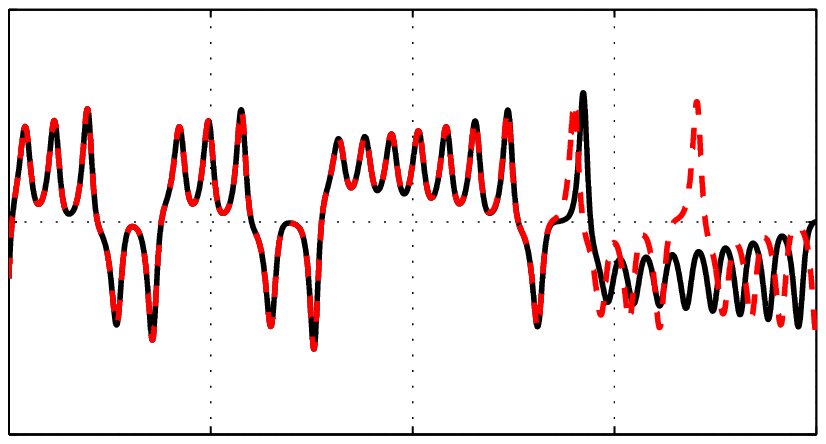}
{\small
% time labels
\put(12.5,2){0}
\put(31.5,2){5}
\put(50,2){10}
\put(69,2){15}
\put(88,2){20}
% y-axis label
\put(7,5){-30}
\put(10,25){0}
\put(8,45){30}
}
\put(4,30){$x$}
\put(45,50){$\eta=0.01$}
\end{overpic}&&
\begin{overpic}[width=0.5\textwidth]{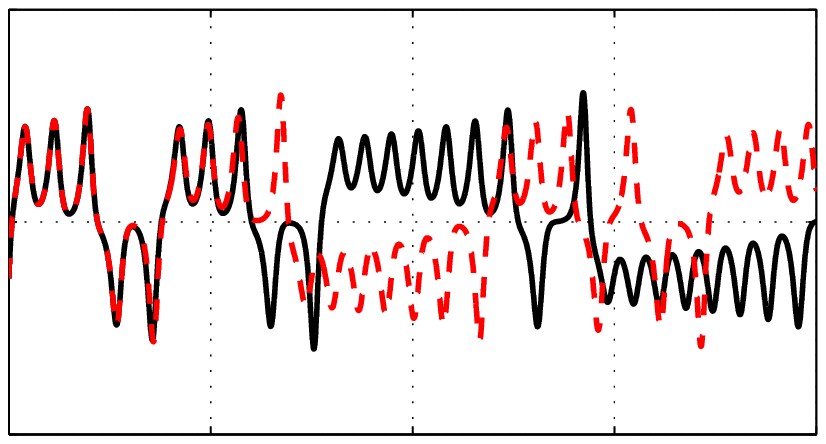}
{\small
% time labels
\put(12.5,2){0}
\put(31.5,2){5}
\put(50,2){10}
\put(69,2){15}
\put(88,2){20}
% y-axis label
\put(7,5){-30}
\put(10,25){0}
\put(8,45){30}
}
\put(4,30){$x$}
\put(46.5,50){$\eta=10$}
\end{overpic}\\
\begin{overpic}[width=0.5\textwidth]{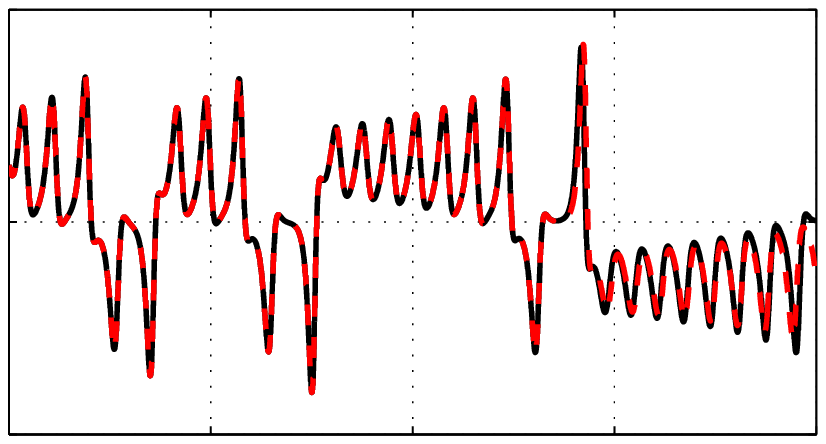}
{\small
% time labels
\put(12.5,2){0}
\put(31.5,2){5}
\put(50,2){10}
\put(69,2){15}
\put(88,2){20}
% y-axis label
\put(7,5){-30}
\put(10,25){0}
\put(8,45){30}
}
\put(4,30){$y$}
\put(47,-2){Time}
\end{overpic}&&
\begin{overpic}[width=0.5\textwidth]{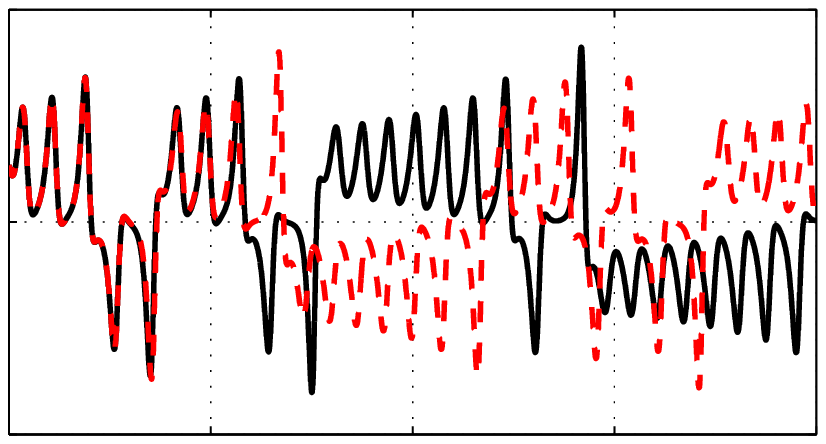}
{\small
% time labels
\put(12.5,2){0}
\put(31.5,2){5}
\put(50,2){10}
\put(69,2){15}
\put(88,2){20}
% y-axis label
\put(7,5){-30}
\put(10,25){0}
\put(8,45){30}
}
\put(4,30){$y$}
\put(47,-2){Time}
\end{overpic}
\end{tabular}
\caption{\small Dynamo view of trajectories of the Lorenz system.  The exact system is shown in black ($-$) and the sparse identified system is shown in the dashed red arrow ({\color{red}$--$}).}\label{Lorenz:dynamo}
\end{center}
\end{figure}

\begin{figure}
\vspace{-.4in}
\begin{center}
\begin{overpic}[width=0.55\textwidth]{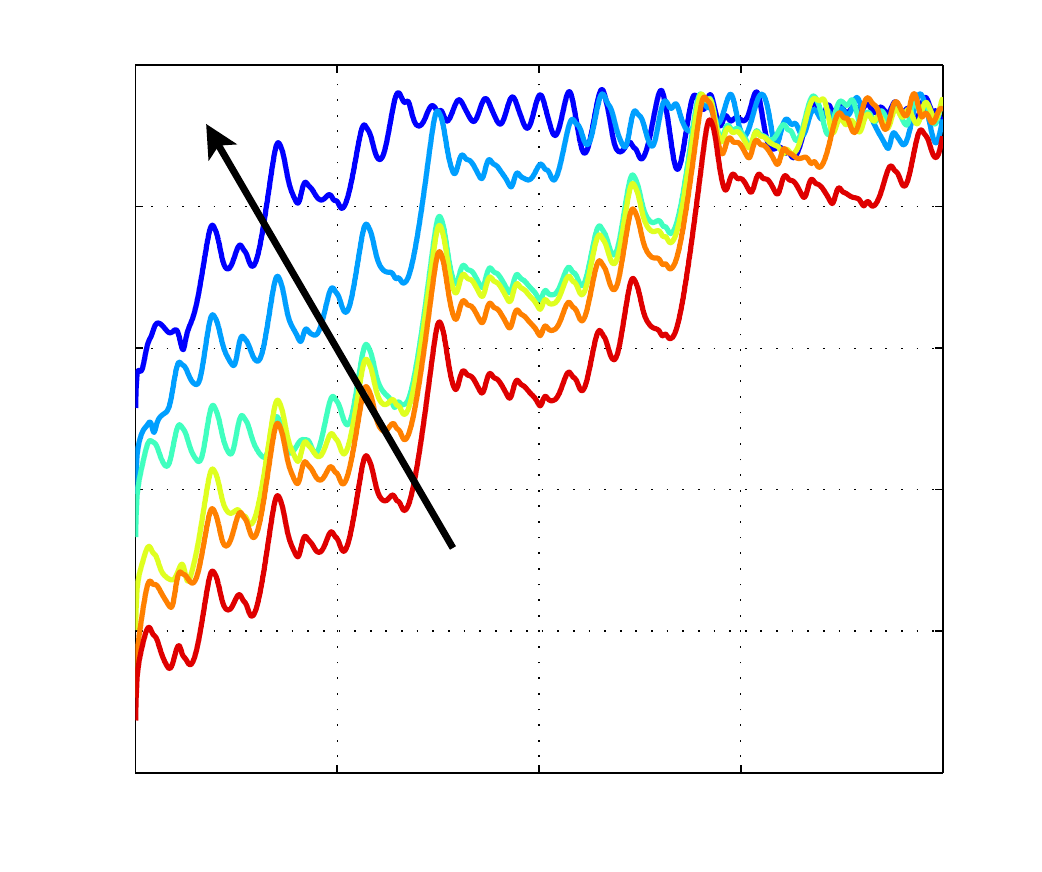}
{\small
\put(12.5,6){0}
\put(31.5,6){5}
\put(50,6){10}
\put(69,6){15}
\put(88,6){20}
\put(4,9){$10^{-8}$}
\put(4,22){$10^{-6}$}
\put(4,36){$10^{-4}$}
\put(4,50){$10^{-2}$}
\put(4,63){$10^0$}
\put(4,75){$10^{2}$}
}
\put(47,1){Time}
\put(-2,38){\begin{sideways}Error\end{sideways}}
\put(45,29){Increasing $\eta$}
\end{overpic}
\caption{\small Error vs. time for sparse identified systems generated from data with increasing sensor noise $\eta$.  This error corresponds to the difference between solid black and dashed red curves in Fig.~\ref{Lorenz:dynamo}.  Sensor noise values are $\eta\in\{0.0001,0.001,0.01,0.1,1.0,10.0\}$.}\label{Lorenz:error}
\end{center}
\end{figure}
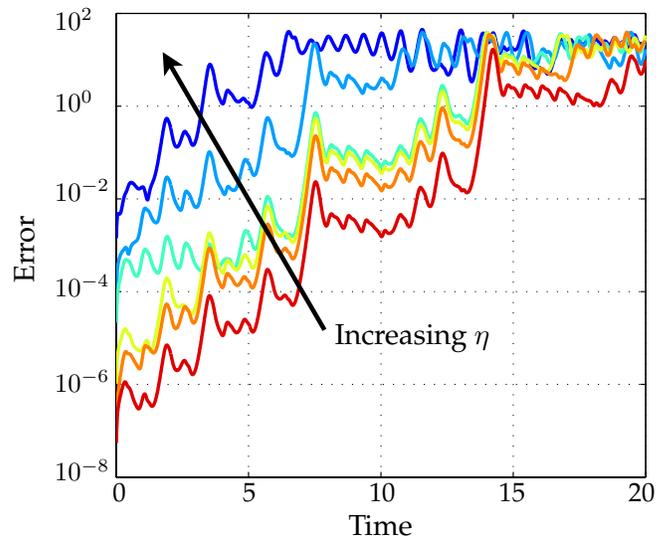

Zero-mean Gaussian measurement noise with variance $\eta$ is added to the derivative calculation to investigate the effect of noisy derivatives.  The short-time ($t=0$ to $t=20$) and long-time ($t=0$ to $t=250$) system reconstruction is shown in Fig.~\ref{Lorenz:Long} for two different noise values, $\eta=0.01$ and $\eta=10$.  The trajectories are also shown in dynamo view in Fig.~\ref{Lorenz:dynamo}, and the $\ell_2$ error vs. time for increasing noise $\eta$ is shown in Fig.~\ref{Lorenz:error}.  Although the $\ell_2$ error increases for large noise values $\eta$, the form of the equations, and hence the attractor dynamics, are accurately captured.  Because the system has a positive Lyapunov exponent, small differences in model coefficients or initial conditions grow exponentially, until saturation, even though the attractor may remain unchanged.

In the Lorenz example, the ability to capture dynamics on the attractor is more important than the ability to predict an individual trajectory, since chaos will quickly cause any small variations in initial conditions or model coefficients to diverge exponentially.  
As shown in Fig.~\ref{FIG00BIG}, our sparse model identification algorithm accurately reproduces the attractor dynamics from chaotic trajectory measurements.  The algorithm not only identifies the correct linear and quadratic terms in the dynamics, but it accurately determines the coefficients to within $.03\%$ of the true values.  When the derivative measurements are contaminated with noise, the correct dynamics are identified, and the attractor is well-preserved for surprisingly large noise values.  When the noise is too large, the structure identification fails before the coefficients become too inaccurate.  
%
%In the Lorenz example, the ability to capture dynamics on the attractor is more important than the ability to predict an individual trajectory, since chaos will quickly cause any small variations in initial conditions or model coefficients to diverge exponentially.  
%As shown in Fig.~\ref{FIG00BIG}, our sparse model identification algorithm correctly identifies the coefficients of the linear and quadratic terms, and accurately reproduces the attractor dynamics from chaotic trajectory measurements.
%

For this example, we use the standard parameters $\sigma=10, \beta=8/3, \rho = 28$, with an initial condition $\begin{bmatrix} x & y & z\end{bmatrix}^T = \begin{bmatrix} -8 & 7 & 27 \end{bmatrix}^T$.  Data is collected from $t=0$ to $t=100$ with a time-step of $\Delta t=0.001$.

\newpage

\subsection{Example 3:  Fluid wake behind a cylinder (Nonlinear PDE)}
% END BIG FIGURE
\begin{figure}[b!]
\begin{center}
\begin{overpic}[width=\textwidth]{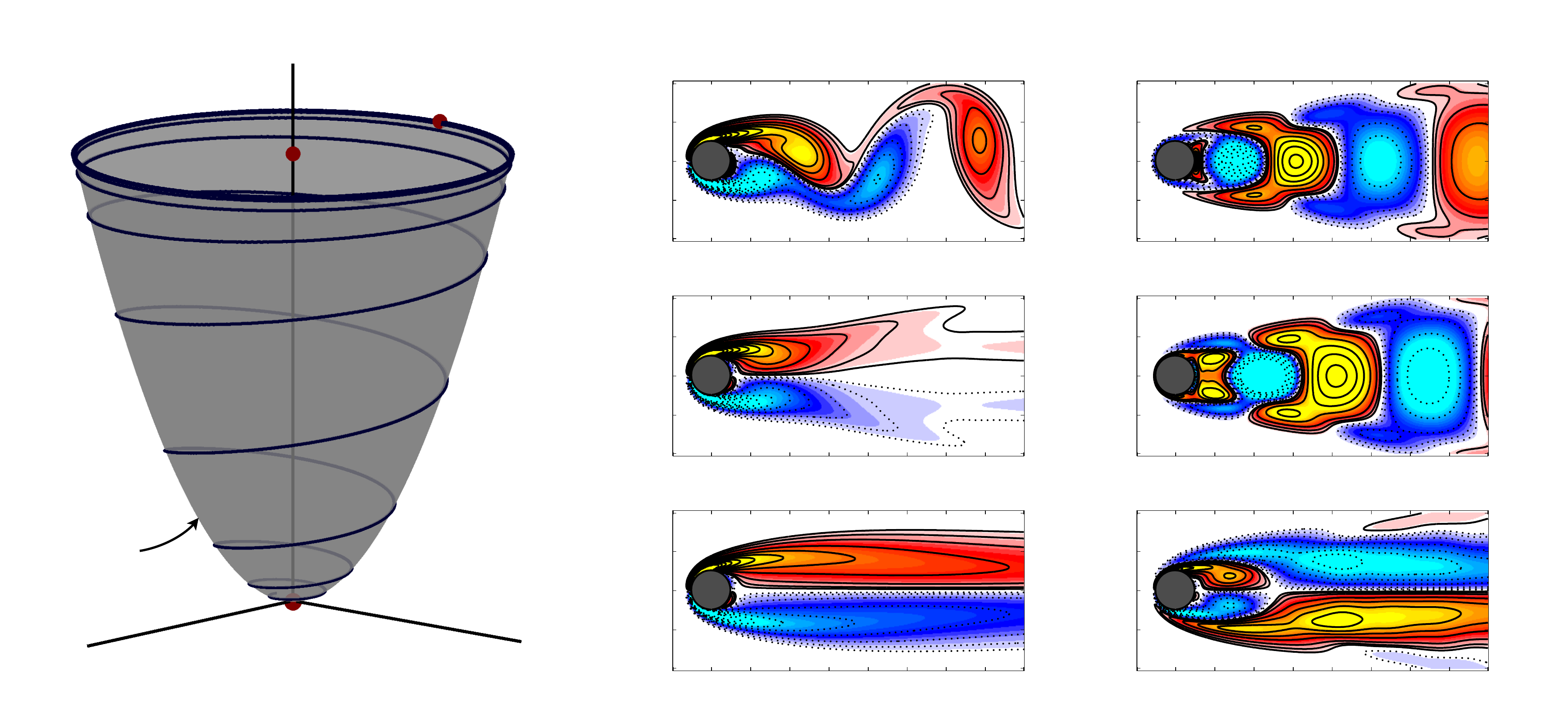}
\put(4,3){$x$}
\put(33.75,3.25){$y$}
\put(19,41.5){$z$}
\put(19,4.5){$C$}
\put(29,38){$A$}
\put(20,34){$B$}

\put(5,39){Limit cycle}
\put(43,40.75){A - vortex shedding}
\put(43,27){B - mean flow}
\put(43,13.25){C - unstable fixed point}
\put(72.5,40.75){$u_x$ - POD mode 1}
\put(72.5,27){$u_y$ - POD mode 2}
\put(72.5,13.25){$u_z$ - shift mode}

\put(3.5,9.5){Slow}
\put(3.5,7.5){manifold}

{\small
\put(42.,0.){-1}
\put(45,0.){0}
\put(47.5,0.){1}
\put(50,0.){2}
\put(52.5,0.){3}
\put(55,0.){4}
\put(57.5,0.){5}
\put(60,0.){6}
\put(62.5,0.){7}
\put(65,0.){8}
% y-labels
\put(40.5,1.75){-2}
\put(40.5,4.25){-1}
\put(41.25,6.75){0}
\put(41.25,9.){1}
\put(41.25,11.5){2}
}
\end{overpic}
\vspace{-.15in}
\caption{\small Illustration of the low-rank dynamics underlying the periodic vortex shedding behind a circular cylinder at low Reynolds number, $\text{Re}=100$.}\label{Cylinder:Manifold}
\vspace{-.3in}
\end{center}
\end{figure}
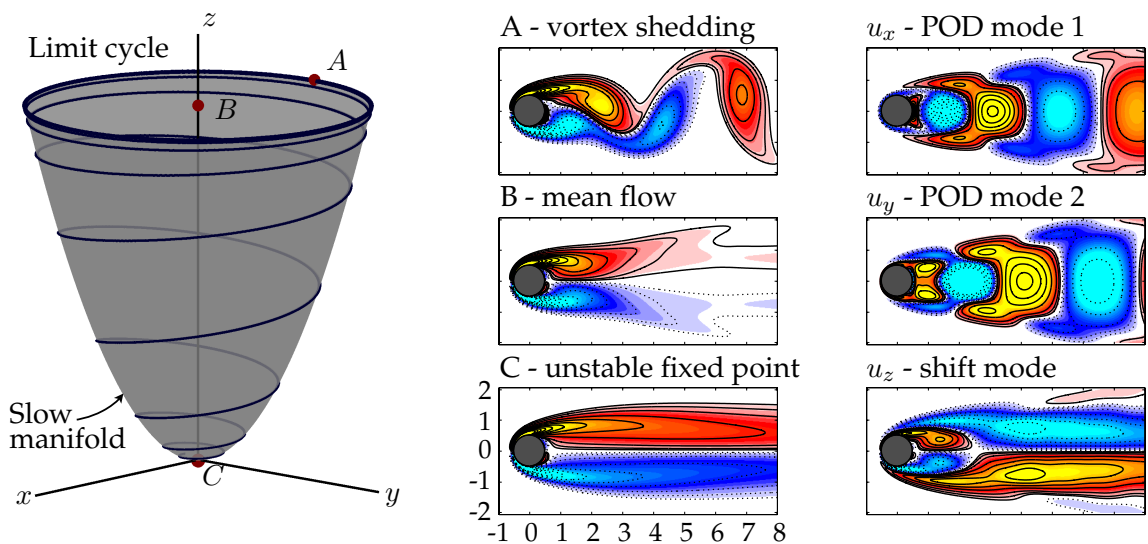
% END BIG FIGURE

The Lorenz system is a low-dimensional model of more realistic high-dimensional partial differential equation (PDE) models for fluid convection in the atmosphere.  
Many systems of interest are governed by PDEs~\cite{Schaeffer2013pnas}, such as weather and climate, epidemiology, and the power grid, to name a few.  
Each of these examples are characterized by big data, consisting of large spatially resolved measurements consisting of millions or billions of states and spanning orders of magnitude of scale in both space and time.  
However, many high-dimensional, real-world systems evolve on a low-dimensional attractor, making the effective dimension much smaller~\cite{HLBR_turb}.  

Here we generalize the sparse identification of nonlinear dynamics method to an example in fluid dynamics that typifies many of the challenges outlined above. 
Data is collected for the fluid flow past a cylinder at Reynolds number 100 using direct numerical simulations of the two-dimensional Navier-Stokes equations~\cite{taira:07ibfs,taira:fastIBPM}.  
Then, the nonlinear dynamic relationship between the dominant coherent structures is identified from these flow field measurements with no knowledge of the governing equations.

The low-Reynolds number flow past a cylinder is a particularly interesting example because of its rich history in fluid mechanics and dynamical systems.  It has long been theorized that turbulence may be the result of a sequence of Hopf bifurcations that occur as the Reynolds number of the flow increases~\cite{Ruelle1971cmp}.  The Reynolds number is a rough measure of the ratio of inertial and viscous forces, and an increasing Reynolds number may correspond, for example, to increasing flow velocity, giving rise to more rich and intricate structures in the fluid.  

After 15 years, the first Hopf bifurcation was discovered in a fluid system, in the transition from a steady laminar wake to laminar periodic vortex shedding at Reynolds number $47$~\cite{Jackson1987jfm,Zebib1987jem,Olinger1988prl}. 
This discovery led to a long-standing debate about how a Hopf bifurcation, with cubic nonlinearity, can be exhibited in a Navier-Stokes fluid with quadratic nonlinearities.  
After 15 more years, this was finally resolved using a separation of time-scales argument and a mean-field model~\cite{Noack2003jfm}, shown in Eq.~\eqref{Eq:MeanFieldModel}.  
It was shown that coupling between oscillatory modes and the base flow gives rise to a slow manifold (see Fig.~\ref{Cylinder:Manifold}, left), which results in algebraic terms that approximate cubic nonlinearities on slow timescales.  

This example provides a compelling test-case for the proposed algorithm, since the underlying form of the dynamics took nearly three decades to uncover.  Indeed, the sparse dynamics algorithm correctly identifies the on-attractor and off-attractor dynamics using quadratic nonlinearities and preserves the correct slow-manifold dynamics.  It is interesting to note that when the off-attractor trajectories are not included in the system identification, the algorithm incorrectly identifies the dynamics using cubic nonlinearities, and fails to correctly identify the dynamics associated with the shift mode, which connects the mean flow to the unstable steady state.  

\subsubsection{Direct numerical simulation}
 The direct numerical simulation involves a fast multi-domain immersed boundary projection method~\cite{taira:07ibfs,taira:fastIBPM}.  Four grids are used, each with a resolution of $450\times 200$, with the finest grid having dimensions of $9\times 4$ cylinder diameters and the largest grid having dimensions of $72\times 32$ diameters.  The finest grid has 90,000 points, and each subsequent coarser grid has 67,500 distinct points.  Thus, if the state includes the vorticity at each grid point, then the state dimension is 292,500.  The vorticity field on the finest grid is shown in Fig.~\ref{Cylinder:Manifold}.  The code is non-dimensionalized so that the cylinder diameter and free-stream velocity are both equal to one:  $D=1$ and $U_{\infty}=1$, respectively.  The simulation time-step is $\Delta t=0.02$ non dimensional time units.

\subsubsection{Mean field model}
To develop a mean-field model for the cylinder wake, first we must reduce the dimension of the system.  The proper orthogonal decomposition (POD)~\cite{HLBR_turb}, provides a low-rank basis that is optimal in the $L^2$ sense, resulting in a hierarchy of orthonormal modes that are ordered by mode energy.  The first two most energetic POD modes capture a significant portion of the energy; the steady-state vortex shedding is a limit cycle in these coordinates.  An additional mode, called the shift mode, is included to capture the transient dynamics connecting the unstable steady state with the mean of the limit cycle~\cite{Noack2003jfm} (i.e., the direction connecting point `C' to point `B' in Fig.~\ref{Cylinder:Manifold}).  

In the three-dimensional coordinate system described above, the mean-field model for the cylinder dynamics are given by:
\begin{subequations}
\label{Eq:MeanFieldModel}
\begin{eqnarray}
\dot{x} & = & \mu x - \omega y + A xz\label{Eq:MeanFieldModel1}\\
\dot{y} & = & \omega x + \mu y + A yz\label{Eq:MeanFieldModel2}\\
\dot{z} & = & -\lambda(z-x^2-y^2).\label{Eq:MeanFieldModel3}
\end{eqnarray}
\end{subequations}
If $\lambda$ is large, so that the $z$-dynamics are fast, then the mean flow rapidly corrects to be on the (slow) manifold $z=x^2+y^2$ given by the amplitude of vortex shedding.  When substituting this algebraic relationship into Eqs.~\ref{Eq:MeanFieldModel1} and~\ref{Eq:MeanFieldModel2}, we recover the Hopf normal form on the slow manifold.

Remarkably, similar dynamics are discovered by the sparse dynamics algorithm, purely from data collected from simulations.  The identified model coefficients, shown in Table~\ref{Tab:Cylinder}, only include quadratic nonlinearities, consistent with the Navier-Stokes equations.  Moreover, the transient behavior, shown in Figs.~\ref{Cylinder:OffAttractor} and~\ref{Cylinder:FarOffAttractor}, is captured qualitatively for solutions that do not start on the slow manifold.  When the off-attractor dynamics in Fig.~\ref{Cylinder:OffAttractor} are not included in the training data, the model incorrectly identifies a simple Hopf normal form in $x$ and $y$ with cubic nonlinearities.

\begin{figure}[where!]
\begin{center}
\begin{tabular}{cccc}
\begin{overpic}[width=.4\textwidth]{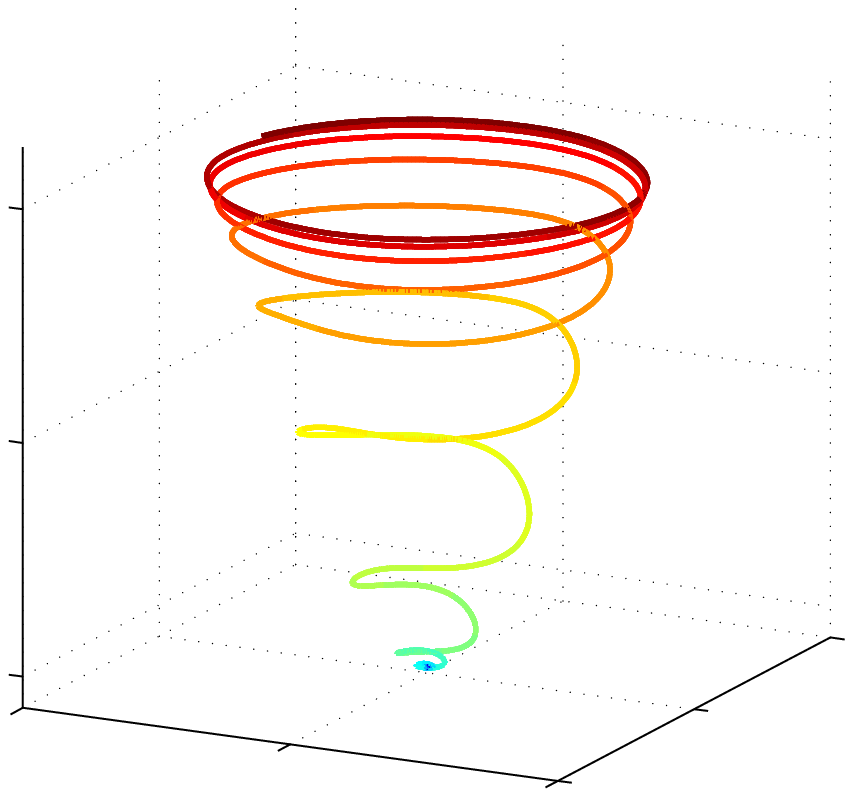}
{\small
%x-label
\put(5,12.5){-200}
\put(35,9.0){0}
\put(58.5,5.5){200}
%y-label
\put(66.5,8.5){-200}
\put(81,15){0}
\put(92.5,22.){200}
%z-label
%\put(5.,19){-50}
%\put(8.25,41.75){0}
%\put(6.25,64.){50}
\put(3.,19){-150}
\put(4.75,41.75){-75}
\put(8.25,64.){0}
}
\put(41,4.5){$x$}
\put(85.5,10.5){$y$}
\put(2.4,52){$z$}
\put(35,87){Full Simulation}
\end{overpic} &&
\begin{overpic}[width=.4\textwidth]{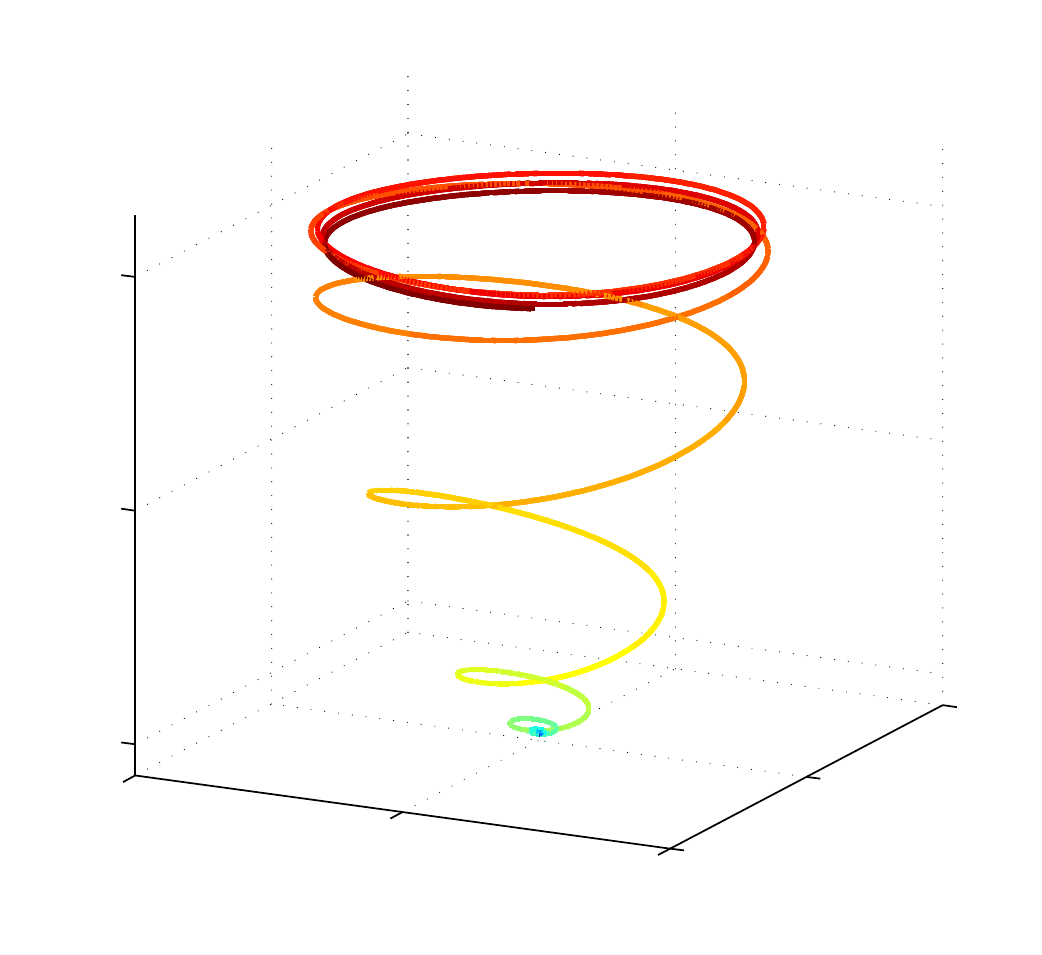}
{\small
%x-label
\put(5,12.5){-200}
\put(35,9.0){0}
\put(58.5,5.5){200}
%y-label
\put(66.5,8.5){-200}
\put(81,15){0}
\put(92.5,22.){200}
%z-label
%\put(5.,19){-50}
%\put(8.25,41.75){0}
%\put(6.25,64.){50}
\put(3.,19){-150}
\put(4.75,41.75){-75}
\put(8.25,64.){0}
}
\put(41,4.5){$x$}
\put(85.5,10.5){$y$}
\put(2.4,52){$z$}
\put(35,87){Identified System}
\end{overpic}
\end{tabular}
\vspace{-.2in}
\caption{\small Evolution of the cylinder wake trajectory in reduced coordinates.  The full simulation (left) comes from direct numerical simulation of the Navier-Stokes equations, and the identified system (right) captures the dynamics on the slow manifold. Color indicates simulation time.}\label{Cylinder:Attractor}
\end{center}
\vspace{-.1in}
\end{figure}

\begin{figure}[where!]
\begin{center}
\begin{tabular}{cccc}
\begin{overpic}[width=.4\textwidth]{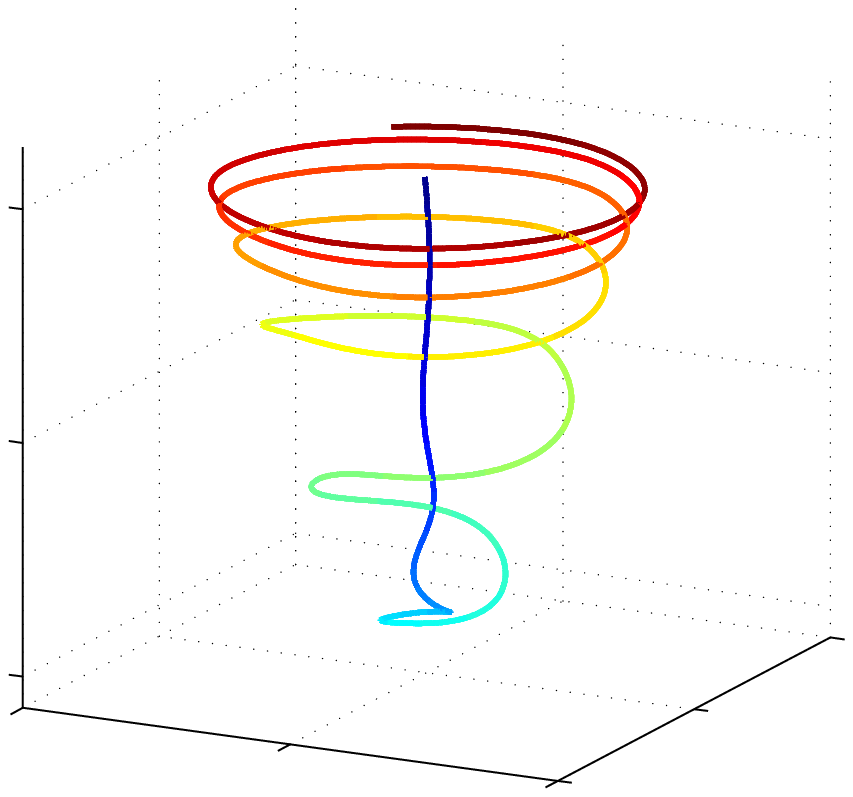}
{\small
%x-label
\put(5,12.5){-200}
\put(35,9.0){0}
\put(58.5,5.5){200}
%y-label
\put(66.5,8.5){-200}
\put(81,15){0}
\put(92.5,22.){200}
%z-label
%\put(5.,19){-50}
%\put(8.25,41.75){0}
%\put(6.25,64.){50}
\put(3.,19){-150}
\put(4.75,41.75){-75}
\put(8.25,64.){0}
}
\put(41,4.5){$x$}
\put(85.5,10.5){$y$}
\put(2.4,52){$z$}
\put(35,87){Full Simulation}
\end{overpic} &&
\begin{overpic}[width=.4\textwidth]{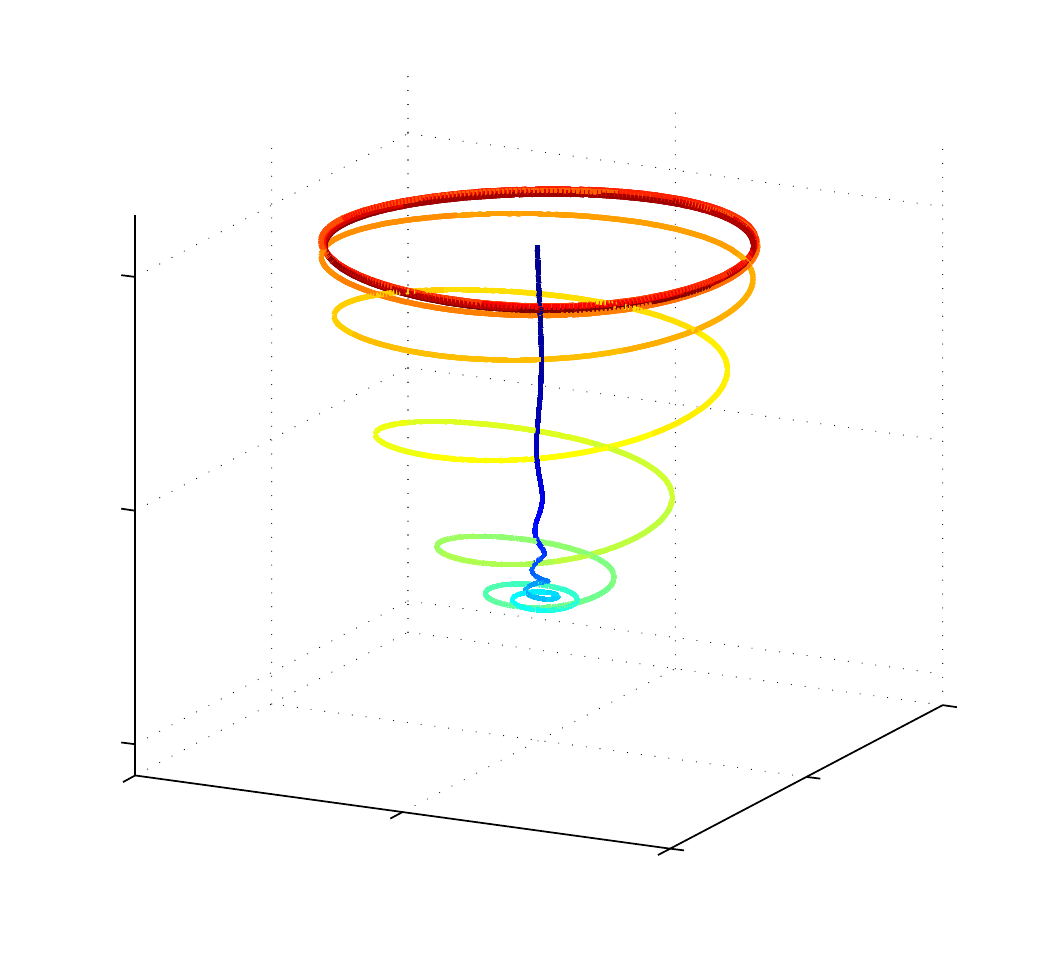}
{\small
%x-label
\put(5,12.5){-200}
\put(35,9.0){0}
\put(58.5,5.5){200}
%y-label
\put(66.5,8.5){-200}
\put(81,15){0}
\put(92.5,22.){200}
%z-label
%\put(5.,19){-50}
%\put(8.25,41.75){0}
%\put(6.25,64.){50}
\put(3.,19){-150}
\put(4.75,41.75){-75}
\put(8.25,64.){0}
}
\put(41,4.5){$x$}
\put(85.5,10.5){$y$}
\put(2.4,52){$z$}
\put(35,87){Identified System}
\end{overpic}
\end{tabular}
\vspace{-.2in}
\caption{\small Evolution of the cylinder wake trajectory starting from a flow state initialized at the mean of the steady-state limit cycle.  Both the full simulation and sparse model capture the off-attractor dynamics, characterized by rapid attraction of the trajectory onto the slow manifold.}\label{Cylinder:OffAttractor}
\end{center}
\vspace{-.1in}
\end{figure}

\begin{figure}[where!]
\begin{center}
\begin{tabular}{cccc}
\begin{overpic}[width=.4\textwidth]{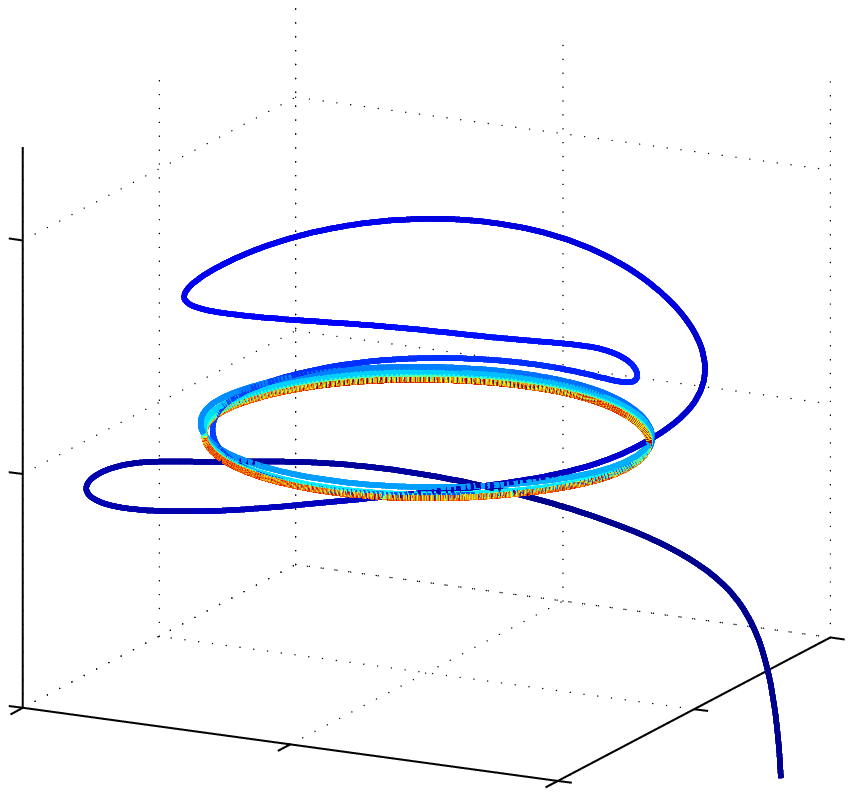}
{\small
%x-label
\put(5,12.5){-200}
\put(35,9.0){0}
\put(58.5,5.5){200}
%y-label
\put(66.5,8.5){-200}
\put(81,15){0}
\put(92.5,22.){200}
%z-label
\put(5.,18){-50}
\put(8.25,38.75){0}
\put(6.25,62.){50}
%\put(3.,19){-150}
%\put(4.75,41.75){-75}
%\put(8.25,64.){0}
}
\put(41,4.5){$x$}
\put(86.5,10.5){$y$}
\put(2.4,52){$z$}
\put(35,87){Full Simulation}
\end{overpic} &&
\begin{overpic}[width=.4\textwidth]{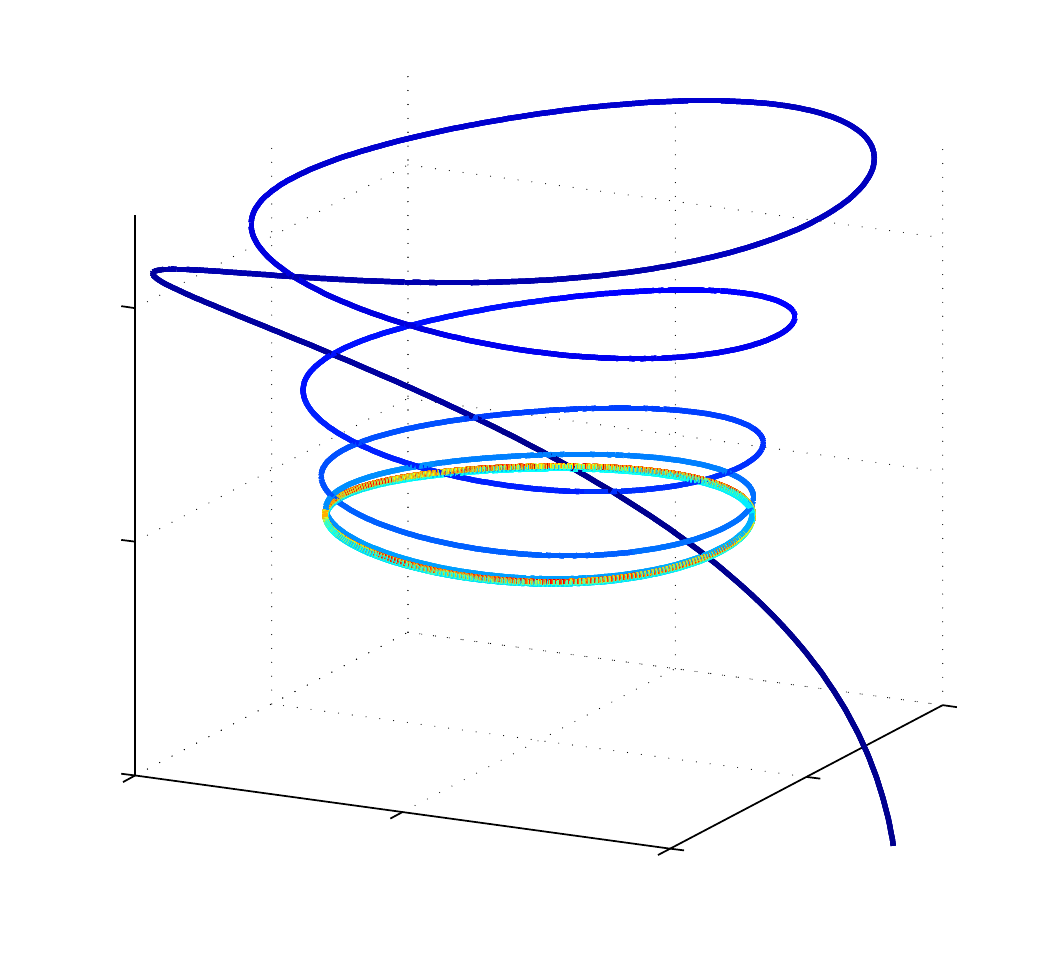}
{\small
%x-label
\put(5,12.5){-200}
\put(35,9.0){0}
\put(58.5,5.5){200}
%y-label
\put(66.5,8.5){-200}
\put(81,15){0}
\put(92.5,22.){200}
%z-label
\put(5.,18){-50}
\put(8.25,38.75){0}
\put(6.25,62.){50}
%\put(3.,19){-150}
%\put(4.75,41.75){-75}
%\put(8.25,64.){0}
}
\put(41,4.5){$x$}
\put(86.5,10.5){$y$}
\put(2.4,52){$z$}
\put(35,87){Identified System}
\end{overpic}
\end{tabular}
\vspace{-.1in}
\caption{\small This simulation corresponds to an initial condition obtained by doubling the magnitude of the limit cycle behavior.  This data was not included in the training of the sparse model.}\label{Cylinder:FarOffAttractor}
\vspace{-.3in}
\end{center}
\vspace{-.3in}
\end{figure}

The data from Fig.~\ref{Cylinder:FarOffAttractor} was not included in the training data, and although qualitatively similar, the identified model does not exactly reproduce the transients.  Since this initial condition had twice the fluctuation energy in the $x$ and $y$ directions, the slow manifold approximation may not be valid here.  Relaxing the sparsity condition, it is possible to obtain models that agree almost perfectly with the data in Figs.~\ref{Cylinder:Attractor}-\ref{Cylinder:FarOffAttractor}, although the model includes higher order nonlinearities.

\subsection{Example 4: Bifurcations and Normal Forms}
It is also possible to identify normal forms associated with a bifurcation parameter $\mu$ by suspending it in the dynamics as a variable:
\begin{subequations}
\begin{eqnarray}
\dot{\mathbf{x}} & = &  \mathbf{f}(\mathbf{x};\mu)\\
\dot\mu & = & 0.
\end{eqnarray}
\end{subequations}
%
%Here we consider the bifurcation parameter $\mu$ as a variable.  
It is then possible to identify the right hand side $\mathbf{f}(\mathbf{x};\mu)$ as a sparse combination of functions of components in $\mathbf{x}$ as well as the bifurcation parameter $\mu$.  This idea is illustrated on two examples, the one-dimensional logistic map and the two-dimensional Hopf normal form.

\subsubsection{Logistic map}
\begin{figure}[where!]
\begin{center}
\begin{tabular}{ccccc}
Stochastic System &&&&  Sparse Identified System\\
\begin{overpic}[width=.3\textwidth]{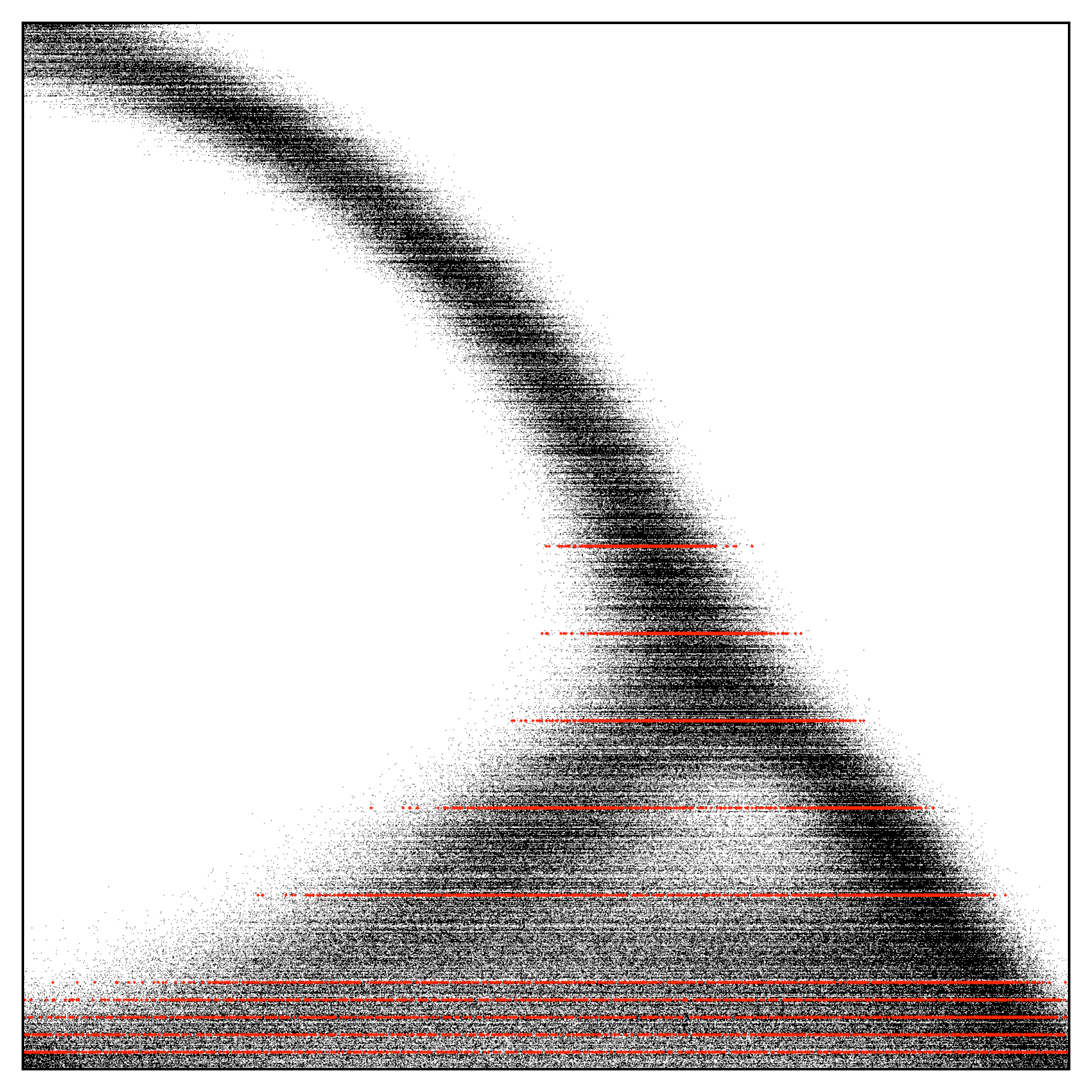}
{\small
\put(1,-4){0}
\put(48.5,-4){0.5}
\put(98,-4){1}
% r-axis
\put(-3,1){4}
\put(-3,33){3}
\put(-3,65){2}
\put(-3,96){1}
}
\put(49.5,-10){$x$}
\put(-10,50){$\mu$}
\end{overpic} && &&
\begin{overpic}[width=.3\textwidth]{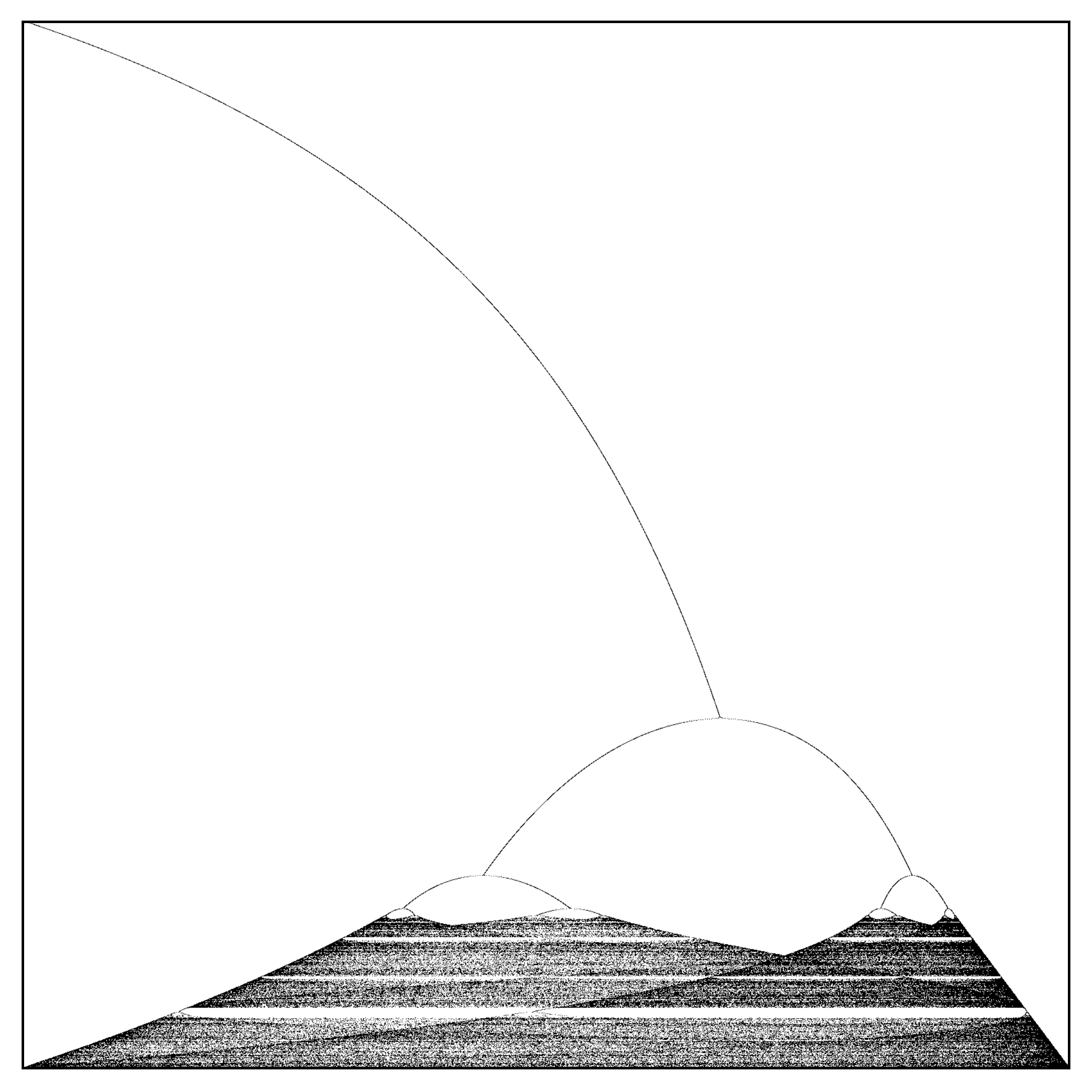}
{\small
\put(1,-4){0}
\put(48.5,-4){0.5}
\put(98,-4){1}
% r-axis
\put(-3,1){4}
\put(-3,33){3}
\put(-3,65){2}
\put(-3,96){1}
}
\put(49.5,-10){$x$}
\put(-12,50){$\mu$}
\end{overpic}
\\
&&&&\\
&&&&\\
%&&&&\\
\begin{overpic}[width=.3\textwidth]{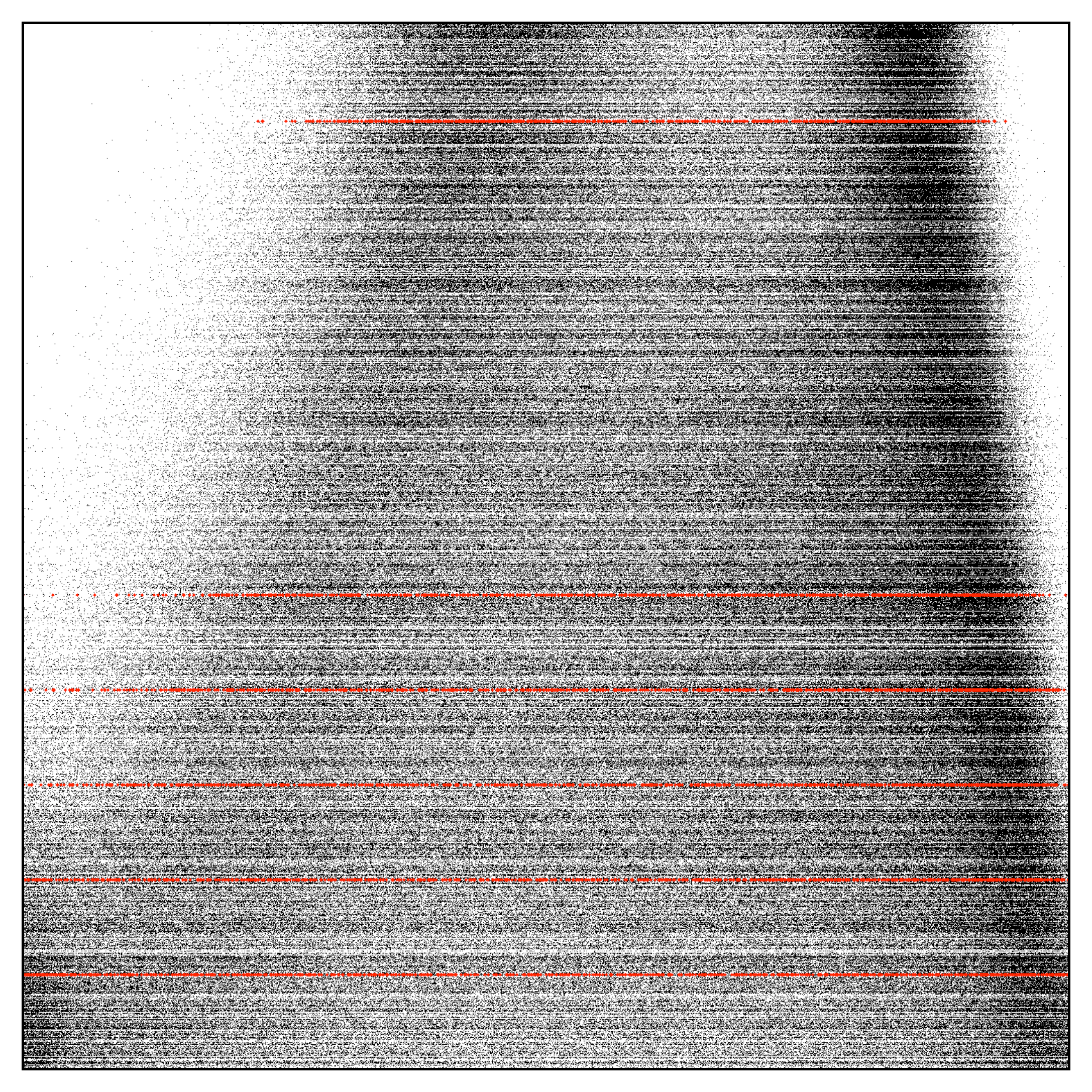}
{\small
\put(1,-4){0}
\put(48.5,-4){0.5}
\put(98,-4){1}
% r-axis
\put(-3,1){4}
\put(-10.5,33){3.82}
\put(-10.5,65){3.63}
\put(-10.5,96){3.45}
}
\put(49.5,-10){$x$}
\put(-12,50){$\mu$}
\end{overpic} &&&&
\begin{overpic}[width=.3\textwidth]{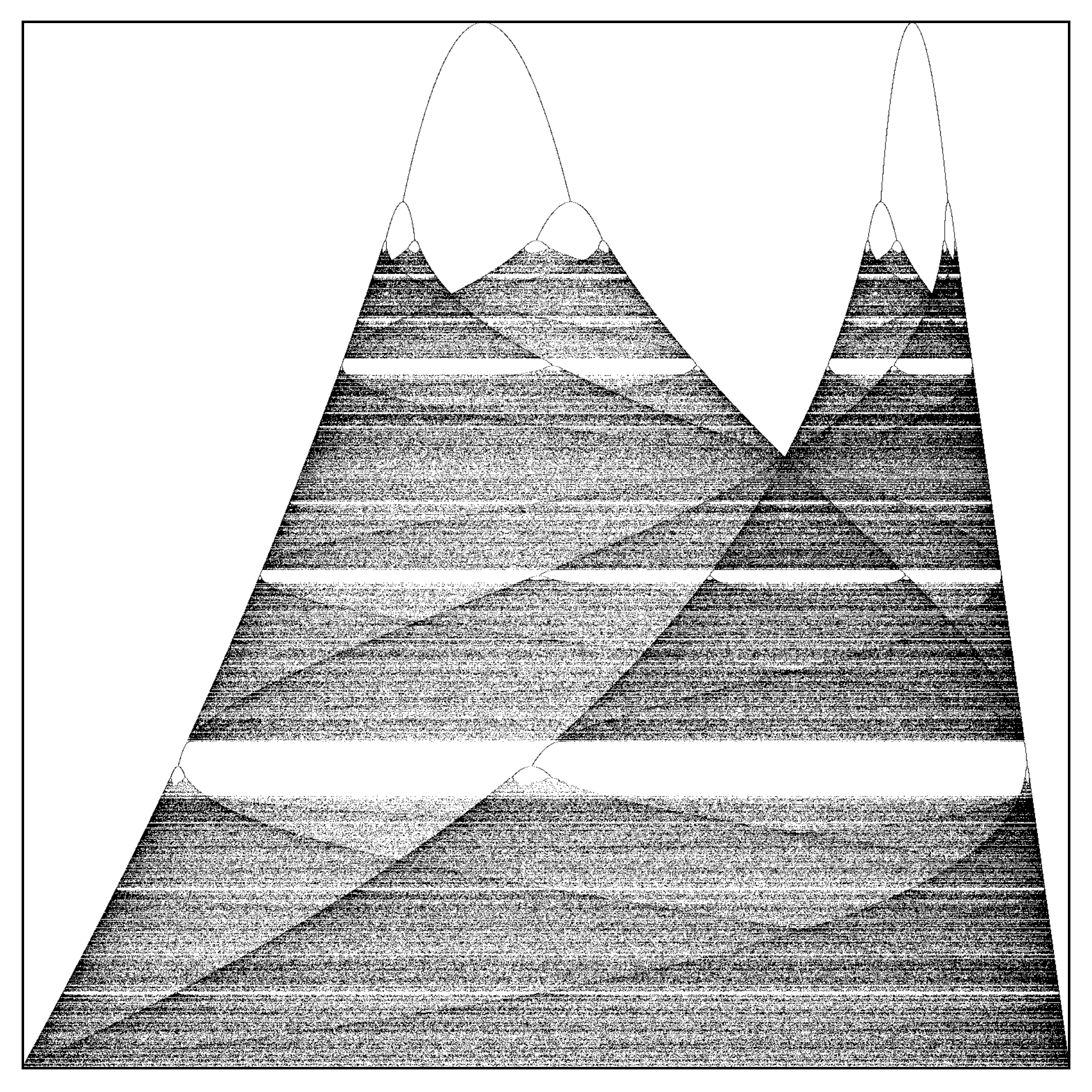}
{\small
\put(1,-4){0}
\put(48.5,-4){0.5}
\put(98,-4){1}
% r-axis
\put(-3,1){4}
\put(-10.5,33){3.82}
\put(-10.5,65){3.63}
\put(-10.5,96){3.45}
}
\put(49.5,-10){$x$}
\put(-10,50){$\mu$}
\end{overpic} \\
&&&&\\
\end{tabular}
\vspace{-.1in}
\caption{\small Attracting sets of the logistic map vs. the parameter $\mu$.  (left) Data from stochastically forced system and (right) the sparse identified system.  Data is sampled at rows indicated in red for $\mu\in\{2.5, 2.75, 3, 3.25, 3.5, 3.75, 3.8, 3.85, 3.9, 3.95\}$.  The forcing $\eta_k$ is Gaussian with magnitude $0.025$.}\label{Fig:Logistic}\vspace{-0.35in}
\end{center}
\end{figure}
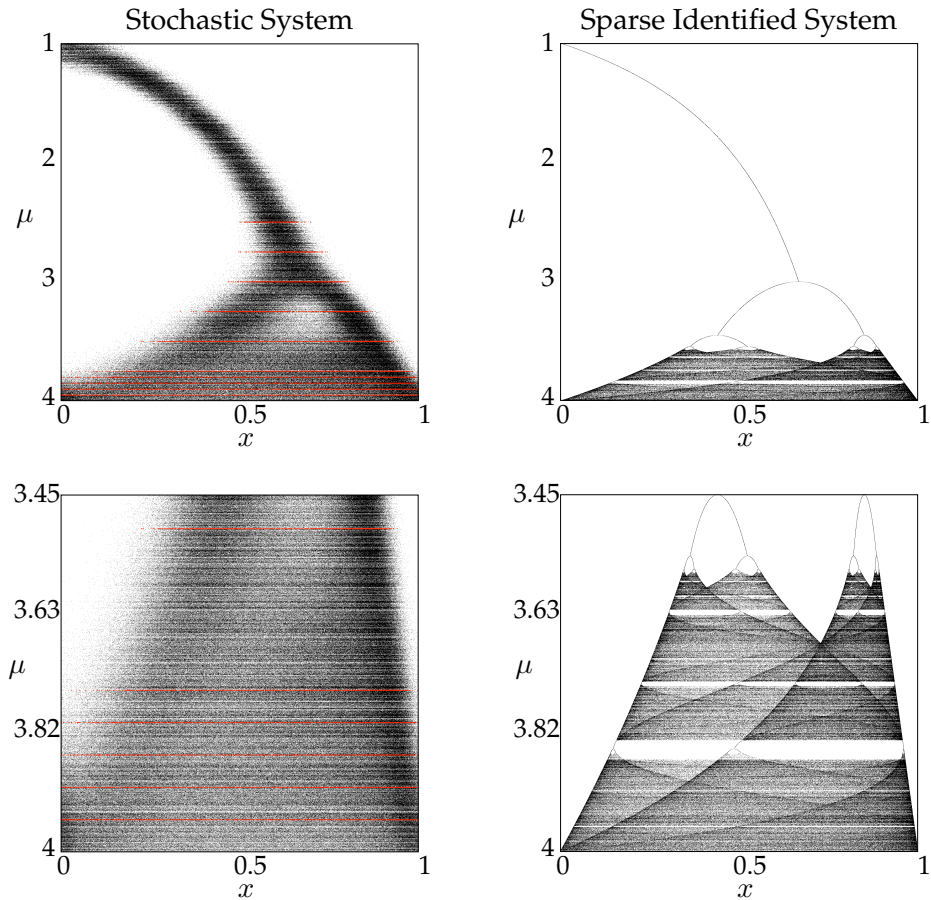

The logistic map is a classical model that exhibits a cascade of bifurcations, leading to chaotic trajectories.  The dynamics with stochastic forcing $\eta_k$ and parameter $\mu$ are given by 
\begin{eqnarray}
x_{k+1} = \mu x_k(1-x_k) + \eta_k.\label{Eq:logistic}
\end{eqnarray}
Sampling the stochastic system at ten parameter values of $\mu$, the algorithm correctly identifies the underlying parameterized dynamics, shown in Fig.~\ref{Fig:Logistic} and Table~\ref{Tab:Logistic}.

\newpage
\subsubsection{Hopf normal form}
\begin{figure}[where!]
\begin{center}
\vspace{-.15in}
\begin{overpic}[width=.75\textwidth]{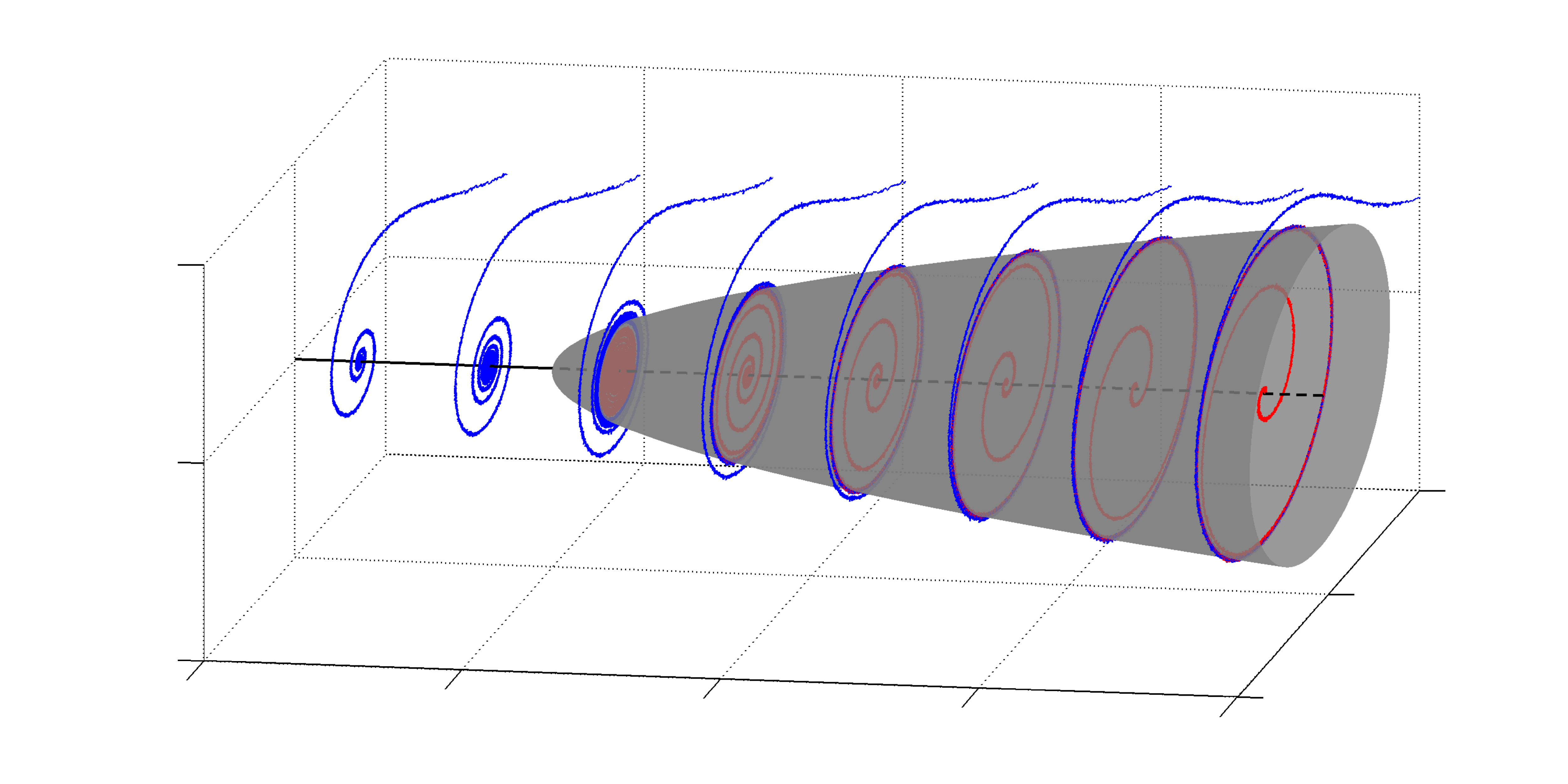}
{\small
% mu axis
\put(9.5,4.75){-0.2}
\put(27.,4.25){0.}
\put(42.5,3.5){0.2}
\put(58.75,2.75){0.4}
\put(75.25,2.25){0.6}
% y axis
\put(8.75,7.5){-1}
\put(9.25,20.15){0}
\put(9.25,32.5){1}
% x axis
\put(81.5,4.75){-1}
\put(87.25,11.5){0}
\put(93,18.){1}
}
\put(35,2.){$\mu$}
\put(6.5,24){$y$}
\put(87,8.){$x$}
\end{overpic}
\vspace{-.1in}
\caption{\small Training data to identify Hopf normal form.  Blue trajectories denote solutions that start outside of the fixed point for $\mu<0$ or the limit cycle for $\mu>0$, and red trajectories denote solutions that start inside of the limit cycle.}\label{Fig:HopfTraining}
\end{center}
\end{figure}

The final example illustrating the ability of the sparse dynamics method to identify parameterized normal forms is the Hopf normal form~\cite{Marsden1976book}.  Noisy data is collected from the Hopf system
\begin{subequations}
\begin{eqnarray}
\dot{x} & =& \mu x + \omega y - Ax(x^2+y^2)\\
\dot{y} & =& -\omega x + \mu y - Ay(x^2 + y^2)
\end{eqnarray}
\end{subequations}
for various values of the parameter $\mu$.  Data is collected on the blue and red trajectories in Fig.~\ref{Fig:HopfTraining}, and noise is added to simulate sensor noise.  The total variation derivative~\cite{Chartrand2011isrnam} is used to de-noise the derivative for use in the algorithm.  

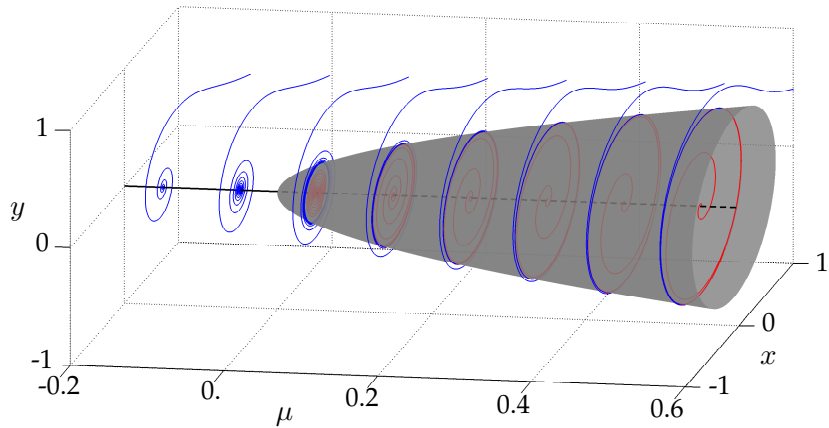
\begin{figure}[where!]
\begin{center}
\vspace{-.25in}
\begin{overpic}[width=.75\textwidth]{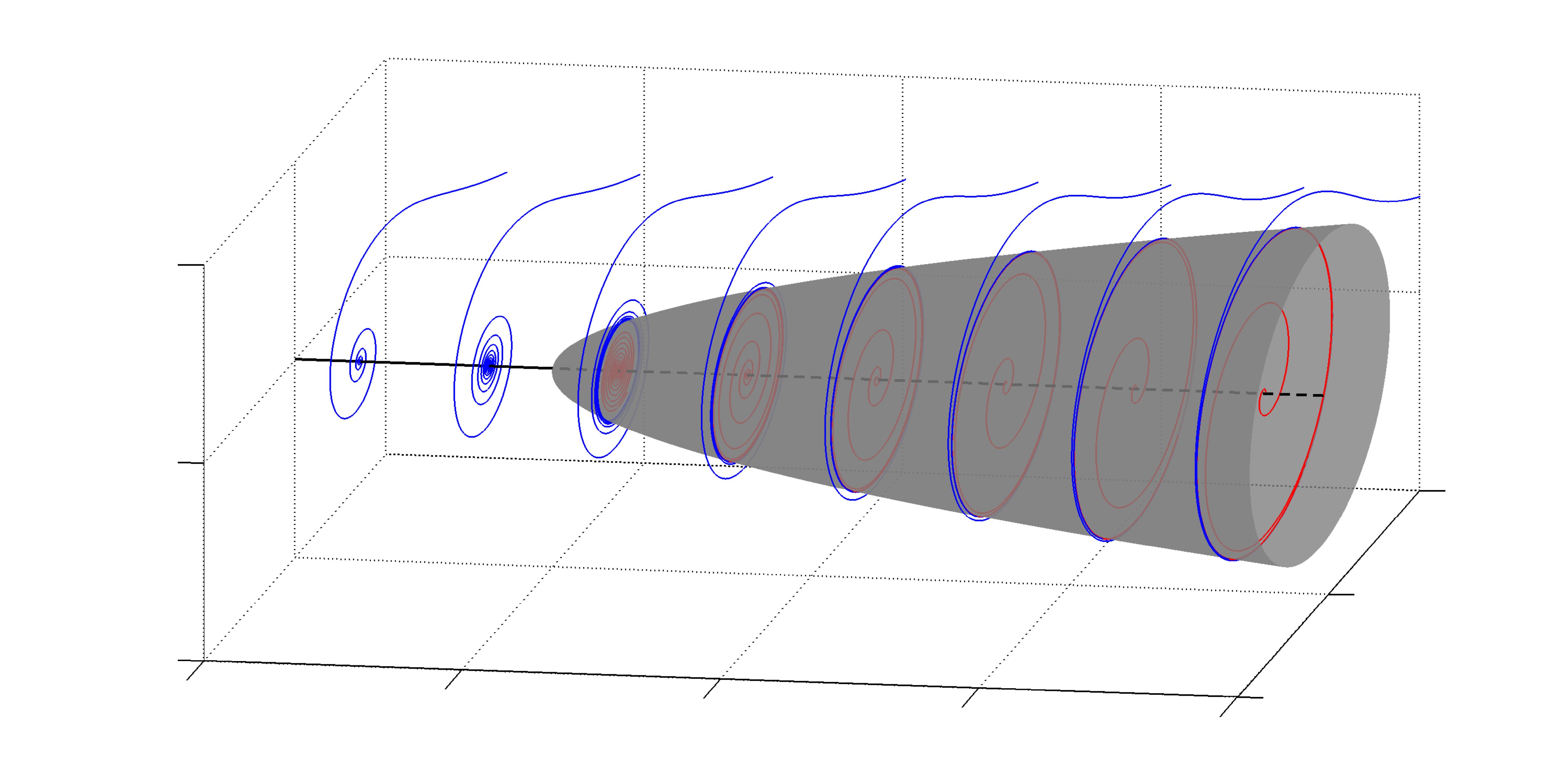}
{\small
% mu axis
\put(9.5,4.75){-0.2}
\put(27.,4.25){0.}
\put(42.5,3.5){0.2}
\put(58.75,2.75){0.4}
\put(75.25,2.25){0.6}
% y axis
\put(8.75,7.5){-1}
\put(9.25,20.15){0}
\put(9.25,32.5){1}
% x axis
\put(81.5,4.75){-1}
\put(87.25,11.5){0}
\put(93,18.){1}
}
\put(35,2.){$\mu$}
\put(6.5,24){$y$}
\put(87,8.){$x$}
\end{overpic}
\vspace{-.1in}
\caption{\small Sparse model captures the Hopf normal form.  Initial conditions are the same as in Fig.~\ref{Fig:HopfTraining}}\label{Fig:HopfReconstruct}
\vspace{-.2in}
\end{center}
\end{figure}

The sparse model identification algorithm correctly identifies the Hopf normal form, with model parameters given in Table~\ref{Tab:Hopf}.  The noise-free model reconstruction is shown in Fig.~\ref{Fig:HopfReconstruct}.  Note that with noise in the training data, although the model terms are correctly identified, the actual values of the cubic terms are off by almost $8\%$.  Collecting more training data or reducing the noise magnitude both improve the model agreement.

%%%%%%%%%%%%
%%% DISCUSSION
%%%%%%%%%%%%

\section{Discussion}
In summary, we have demonstrated a powerful new technique to identify nonlinear dynamical systems from data without assumptions on the form of the governing equations.  
This builds on prior work in symbolic regression but with innovations related to sparse regression, which allow our algorithms to scale to high-dimensional systems.  
We demonstrate this new method on a number of example systems exhibiting chaos, high-dimensional data with low-rank coherence, and parameterized dynamics.  
As shown in the Lorenz example, the ability to predict a specific trajectory may be less important than the ability to capture the attractor dynamics.  
The example from fluid dynamics highlights the remarkable ability of this method to extract dynamics in a fluid system that took three decades for experts in the community to explain.  
There are numerous fields where this method may be applied, where there is ample data and the absence of governing equations, including neuroscience, climate science, epidemiology, and financial markets. 
Fields that already use genetic programming, such as machine learning control for turbulent fluid systems~\cite{Brunton2015amr,Parezanovic2014JFTC}, may also benefit.  
Finally, normal forms may be discovered by including parameters in the optimization, as shown on two examples.  
The identification of \emph{sparse} governing equations and parameterizations marks a significant step toward the long-held goal of intelligent, unassisted identification of dynamical systems.

A number of open problems remain surrounding the dynamical systems aspects of this procedure.  For example, many systems possess dynamical symmetries and conserved quantities that may alter the form of the identified dynamics.  For example, the degenerate identification of a linear system in a space of high-order polynomial nonlinearities suggest a connection with near-identity transformations and dynamic similarity.  We believe that this may be a fruitful line of research.  Finally, it will be important to identify which approximating function space to use based on the data available.  For example, it may be possible to improve the function space to make the dynamics more sparse through subsequent coordinate transformations~\cite{guckenheimer_holmes}.

Data science is not a panacea for all problems in science and engineering, but used in the right way, it provides a principled approach to maximally leverage the data that we have and inform what new data to collect.  
Big data is happening all across the sciences, where the data is inherently \emph{dynamic}, and where traditional approaches are prone to overfitting.  
Data discovery algorithms that produce \emph{parsimonious} models are both rare and desirable.
Data-science will only become more critical to efforts in science in engineering, where data is abundant, but physical laws remain elusive.
These efforts include understanding the neural basis of cognition, extracting and predicting coherent changes in the climate, stabilizing financial markets, managing the spread of disease, and controlling turbulence,

\section*{Acknowledgements}
We gratefully acknowledge valuable discussions with Bingni W. Brunton and Bernd R. Noack.  SLB acknowledges support from the University of Washington department of Mechanical Engineering and as a Data Science Fellow in the eScience Institute (NSF, Moore-Sloan Foundation, Washington Research Foundation).  JLP thanks Bill and Melinda Gates for their active support of the Institute for Disease Modeling and their sponsorship through the Global Good Fund.  JNK acknowledges support from the  U.S. Air Force Office of Scientific Research (FA9550-09-0174).

\newpage
%%%%%%%%%%%%
%%% BIBLIOGRAPHY
%%%%%%%%%%%%
\begin{spacing}{.89}
\footnotesize{
\bibliographystyle{plain}
\bibliography{references}
}
\end{spacing}

\newpage
\section*{Appendix}
%%%% LINEAR DAMPED HARMONIC OSCILLATOR

\begin{table}[where!]
\caption{\small Damped harmonic oscillator with linear terms.}\label{Tab:Ex1_2dLin}
\begin{center}
{\footnotesize
\begin{verbatim}
    ''         'xdot'       'ydot'   
    '1'        [      0]    [      0]
    'x'        [-0.1015]    [-1.9990]
    'y'        [ 2.0027]    [-0.0994]
    'xx'       [      0]    [      0]
    'xy'       [      0]    [      0]
    'yy'       [      0]    [      0]
    'xxx'      [      0]    [      0]
    'xxy'      [      0]    [      0]
    'xyy'      [      0]    [      0]
    'yyy'      [      0]    [      0]
    'xxxx'     [      0]    [      0]
    'xxxy'     [      0]    [      0]
    'xxyy'     [      0]    [      0]
    'xyyy'     [      0]    [      0]
    'yyyy'     [      0]    [      0]
    'xxxxx'    [      0]    [      0]
    'xxxxy'    [      0]    [      0]
    'xxxyy'    [      0]    [      0]
    'xxyyy'    [      0]    [      0]
    'xyyyy'    [      0]    [      0]
    'yyyyy'    [      0]    [      0]
        \end{verbatim}
    }
    \end{center}
    \end{table}

%%%% CUBIC DAMPED HARMONIC OSCILLATOR

\begin{table}[where!]
\caption{\small Damped harmonic oscillator with cubic nonlinearity.}\label{Tab:Ex1_2dCub}
\begin{center}
{\footnotesize
\begin{verbatim}
    ''         'xdot'       'ydot'   
    '1'        [      0]    [      0]
    'x'        [      0]    [      0]
    'y'        [      0]    [      0]
    'xx'       [      0]    [      0]
    'xy'       [      0]    [      0]
    'yy'       [      0]    [      0]
    'xxx'      [-0.0996]    [-1.9994]
    'xxy'      [      0]    [      0]
    'xyy'      [      0]    [      0]
    'yyy'      [ 1.9970]    [-0.0979]
    'xxxx'     [      0]    [      0]
    'xxxy'     [      0]    [      0]
    'xxyy'     [      0]    [      0]
    'xyyy'     [      0]    [      0]
    'yyyy'     [      0]    [      0]
    'xxxxx'    [      0]    [      0]
    'xxxxy'    [      0]    [      0]
    'xxxyy'    [      0]    [      0]
    'xxyyy'    [      0]    [      0]
    'xyyyy'    [      0]    [      0]
    'yyyyy'    [      0]    [      0]

    \end{verbatim}
    }
    \end{center}
    \end{table}

 %%%% THIRD ORDER LINEAR SYSTEM

\begin{table}[where!]
\caption{\small Three-dimensional linear system.}\label{Tab:Ex1_3d}
\begin{center}
{\footnotesize
\begin{verbatim}
    ''      'xdot'       'ydot'       'zdot'   
    '1'     [      0]    [      0]    [      0]
    'x'     [-0.0996]    [-1.9997]    [      0]
    'y'     [ 2.0005]    [-0.0994]    [      0]
    'z'     [      0]    [      0]    [-0.3003]
    'xx'    [      0]    [      0]    [      0]
    'xy'    [      0]    [      0]    [      0]
    'xz'    [      0]    [      0]    [      0]
    'yy'    [      0]    [      0]    [      0]
    'yz'    [      0]    [      0]    [      0]
    'zz'    [      0]    [      0]    [      0]
    \end{verbatim}
    }
    \end{center}
    \end{table}

%%%% LORENZ SYSTEM
\begin{table}[where]
\caption{\small Lorenz system identified using sparse representation with $\eta=1.0$.}
\begin{center}
{\footnotesize
\begin{verbatim}
    ''         'xdot'       'ydot'       'zdot'   
    '1'        [      0]    [      0]    [      0]
    'x'        [-9.9996]    [27.9980]    [      0]
    'y'        [ 9.9998]    [-0.9997]    [      0]
    'z'        [      0]    [      0]    [-2.6665]
    'xx'       [      0]    [      0]    [      0]
    'xy'       [      0]    [      0]    [ 1.0000]
    'xz'       [      0]    [-0.9999]    [      0]
    'yy'       [      0]    [      0]    [      0]
    'yz'       [      0]    [      0]    [      0]
    'zz'       [      0]    [      0]    [      0]
    'xxx'      [      0]    [      0]    [      0]
    'xxy'      [      0]    [      0]    [      0]
    'xxz'      [      0]    [      0]    [      0]
    'xyy'      [      0]    [      0]    [      0]
    'xyz'      [      0]    [      0]    [      0]
    'xzz'      [      0]    [      0]    [      0]
    'yyy'      [      0]    [      0]    [      0]
    'yyz'      [      0]    [      0]    [      0]
    'yzz'      [      0]    [      0]    [      0]
    'zzz'      [      0]    [      0]    [      0]
    'xxxx'     [      0]    [      0]    [      0]
    'xxxy'     [      0]    [      0]    [      0]
    'xxxz'     [      0]    [      0]    [      0]
    'xxyy'     [      0]    [      0]    [      0]
    'xxyz'     [      0]    [      0]    [      0]
    'xxzz'     [      0]    [      0]    [      0]
    'xyyy'     [      0]    [      0]    [      0]
    'xyyz'     [      0]    [      0]    [      0]
    'xyzz'     [      0]    [      0]    [      0]
    'xzzz'     [      0]    [      0]    [      0]
    'yyyy'     [      0]    [      0]    [      0]
    'yyyz'     [      0]    [      0]    [      0]
    'yyzz'     [      0]    [      0]    [      0]
    'yzzz'     [      0]    [      0]    [      0]
    'zzzz'     [      0]    [      0]    [      0]
    'xxxxx'    [      0]    [      0]    [      0]
    'xxxxy'    [      0]    [      0]    [      0]
    'xxxxz'    [      0]    [      0]    [      0]
    'xxxyy'    [      0]    [      0]    [      0]
    'xxxyz'    [      0]    [      0]    [      0]
    'xxxzz'    [      0]    [      0]    [      0]
    'xxyyy'    [      0]    [      0]    [      0]
    'xxyyz'    [      0]    [      0]    [      0]
    'xxyzz'    [      0]    [      0]    [      0]
    'xxzzz'    [      0]    [      0]    [      0]
    'xyyyy'    [      0]    [      0]    [      0]
    'xyyyz'    [      0]    [      0]    [      0]
    'xyyzz'    [      0]    [      0]    [      0]
    'xyzzz'    [      0]    [      0]    [      0]
    'xzzzz'    [      0]    [      0]    [      0]
    'yyyyy'    [      0]    [      0]    [      0]
    'yyyyz'    [      0]    [      0]    [      0]
    'yyyzz'    [      0]    [      0]    [      0]
    'yyzzz'    [      0]    [      0]    [      0]
    'yzzzz'    [      0]    [      0]    [      0]
    'zzzzz'    [      0]    [      0]    [      0]
\end{verbatim}
}
\end{center}
\end{table}

%%%% CYLINDER EXAMPLE

\begin{table}[where!]
\caption{\small Dynamics of cylinder wake modes using sparse representation.  Notice that quadratic terms are identified.}\label{Tab:Cylinder}
\begin{center}
\vspace{-.15in}
{\footnotesize
\begin{verbatim}
    ''         'xdot'           'ydot'           'zdot'       
    '1'        [    -0.1225]    [    -0.0569]    [   -20.8461]
    'x'        [    -0.0092]    [     1.0347]    [-4.6476e-04]
    'y'        [    -1.0224]    [     0.0047]    [ 2.4057e-04]
    'z'        [-9.2203e-04]    [-4.4932e-04]    [    -0.2968]
    'xx'       [          0]    [          0]    [     0.0011]
    'xy'       [          0]    [          0]    [          0]
    'xz'       [ 2.1261e-04]    [     0.0022]    [          0]
    'yy'       [          0]    [          0]    [ 8.6432e-04]
    'yz'       [    -0.0019]    [    -0.0018]    [          0]
    'zz'       [          0]    [          0]    [    -0.0010]
    'xxx'      [          0]    [          0]    [          0]
    'xxy'      [          0]    [          0]    [          0]
    'xxz'      [          0]    [          0]    [          0]
    'xyy'      [          0]    [          0]    [          0]
    'xyz'      [          0]    [          0]    [          0]
    'xzz'      [          0]    [          0]    [          0]
    'yyy'      [          0]    [          0]    [          0]
    'yyz'      [          0]    [          0]    [          0]
    'yzz'      [          0]    [          0]    [          0]
    'zzz'      [          0]    [          0]    [          0]
    'xxxx'     [          0]    [          0]    [          0]
    'xxxy'     [          0]    [          0]    [          0]
    'xxxz'     [          0]    [          0]    [          0]
    'xxyy'     [          0]    [          0]    [          0]
    'xxyz'     [          0]    [          0]    [          0]
    'xxzz'     [          0]    [          0]    [          0]
    'xyyy'     [          0]    [          0]    [          0]
    'xyyz'     [          0]    [          0]    [          0]
    'xyzz'     [          0]    [          0]    [          0]
    'xzzz'     [          0]    [          0]    [          0]
    'yyyy'     [          0]    [          0]    [          0]
    'yyyz'     [          0]    [          0]    [          0]
    'yyzz'     [          0]    [          0]    [          0]
    'yzzz'     [          0]    [          0]    [          0]
    'zzzz'     [          0]    [          0]    [          0]
    'xxxxx'    [          0]    [          0]    [          0]
    'xxxxy'    [          0]    [          0]    [          0]
    'xxxxz'    [          0]    [          0]    [          0]
    'xxxyy'    [          0]    [          0]    [          0]
    'xxxyz'    [          0]    [          0]    [          0]
    'xxxzz'    [          0]    [          0]    [          0]
    'xxyyy'    [          0]    [          0]    [          0]
    'xxyyz'    [          0]    [          0]    [          0]
    'xxyzz'    [          0]    [          0]    [          0]
    'xxzzz'    [          0]    [          0]    [          0]
    'xyyyy'    [          0]    [          0]    [          0]
    'xyyyz'    [          0]    [          0]    [          0]
    'xyyzz'    [          0]    [          0]    [          0]
    'xyzzz'    [          0]    [          0]    [          0]
    'xzzzz'    [          0]    [          0]    [          0]
    'yyyyy'    [          0]    [          0]    [          0]
    'yyyyz'    [          0]    [          0]    [          0]
    'yyyzz'    [          0]    [          0]    [          0]
    'yyzzz'    [          0]    [          0]    [          0]
    'yzzzz'    [          0]    [          0]    [          0]
    'zzzzz'    [          0]    [          0]    [          0]
\end{verbatim}
}
\end{center}
\end{table}

%%%% LOGISTIC NORMAL FORM

\begin{table}[where!]
\caption{\small Logistic map identified using sparse representation.}\label{Tab:Logistic}
\begin{center}
{\footnotesize
\begin{verbatim}
    ''         'x_{k+1}'    'r_{k+1}'
    '1'        [      0]    [     0]
    'x'        [      0]    [     0]
    'r'        [      0]    [1.0000]
    'xx'       [      0]    [     0]
    'xr'       [ 0.9993]    [     0]
    'rr'       [      0]    [     0]
    'xxx'      [      0]    [     0]
    'xxr'      [-0.9989]    [     0]
    'xrr'      [      0]    [     0]
    'rrr'      [      0]    [     0]
    'xxxx'     [      0]    [     0]
    'xxxr'     [      0]    [     0]
    'xxrr'     [      0]    [     0]
    'xrrr'     [      0]    [     0]
    'rrrr'     [      0]    [     0]
    'xxxxx'    [      0]    [     0]
    'xxxxr'    [      0]    [     0]
    'xxxrr'    [      0]    [     0]
    'xxrrr'    [      0]    [     0]
    'xrrrr'    [      0]    [     0]
    'rrrrr'    [      0]    [     0]
    \end{verbatim}
    }
    \end{center}
    \end{table}

%%%% HOPF NORMAL FORM

\begin{table}[where!]
\caption{\small Hopf normal form identified using sparse representation.  Here \texttt{u} represents the bifurcation parameter $\mu$.}\label{Tab:Hopf}
\begin{center}
\vspace{-.15in}
{\footnotesize
\begin{verbatim}
    ''         'xdot'       'ydot'       'udot'
    '1'        [      0]    [      0]    [   0]
    'x'        [      0]    [ 0.9914]    [   0]
    'y'        [-0.9920]    [      0]    [   0]
    'u'        [      0]    [      0]    [   0]
    'xx'       [      0]    [      0]    [   0]
    'xy'       [      0]    [      0]    [   0]
    'xu'       [ 0.9269]    [      0]    [   0]
    'yy'       [      0]    [      0]    [   0]
    'yu'       [      0]    [ 0.9294]    [   0]
    'uu'       [      0]    [      0]    [   0]
    'xxx'      [-0.9208]    [      0]    [   0]
    'xxy'      [      0]    [-0.9244]    [   0]
    'xxu'      [      0]    [      0]    [   0]
    'xyy'      [-0.9211]    [      0]    [   0]
    'xyu'      [      0]    [      0]    [   0]
    'xuu'      [      0]    [      0]    [   0]
    'yyy'      [      0]    [-0.9252]    [   0]
    'yyu'      [      0]    [      0]    [   0]
    'yuu'      [      0]    [      0]    [   0]
    'uuu'      [      0]    [      0]    [   0]
    'xxxx'     [      0]    [      0]    [   0]
    'xxxy'     [      0]    [      0]    [   0]
    'xxxu'     [      0]    [      0]    [   0]
    'xxyy'     [      0]    [      0]    [   0]
    'xxyu'     [      0]    [      0]    [   0]
    'xxuu'     [      0]    [      0]    [   0]
    'xyyy'     [      0]    [      0]    [   0]
    'xyyu'     [      0]    [      0]    [   0]
    'xyuu'     [      0]    [      0]    [   0]
    'xuuu'     [      0]    [      0]    [   0]
    'yyyy'     [      0]    [      0]    [   0]
    'yyyu'     [      0]    [      0]    [   0]
    'yyuu'     [      0]    [      0]    [   0]
    'yuuu'     [      0]    [      0]    [   0]
    'uuuu'     [      0]    [      0]    [   0]
    'xxxxx'    [      0]    [      0]    [   0]
    'xxxxy'    [      0]    [      0]    [   0]
    'xxxxu'    [      0]    [      0]    [   0]
    'xxxyy'    [      0]    [      0]    [   0]
    'xxxyu'    [      0]    [      0]    [   0]
    'xxxuu'    [      0]    [      0]    [   0]
    'xxyyy'    [      0]    [      0]    [   0]
    'xxyyu'    [      0]    [      0]    [   0]
    'xxyuu'    [      0]    [      0]    [   0]
    'xxuuu'    [      0]    [      0]    [   0]
    'xyyyy'    [      0]    [      0]    [   0]
    'xyyyu'    [      0]    [      0]    [   0]
    'xyyuu'    [      0]    [      0]    [   0]
    'xyuuu'    [      0]    [      0]    [   0]
    'xuuuu'    [      0]    [      0]    [   0]
    'yyyyy'    [      0]    [      0]    [   0]
    'yyyyu'    [      0]    [      0]    [   0]
    'yyyuu'    [      0]    [      0]    [   0]
    'yyuuu'    [      0]    [      0]    [   0]
    'yuuuu'    [      0]    [      0]    [   0]
    'uuuuu'    [      0]    [      0]    [   0]
\end{verbatim}
}
\end{center}
\end{table}

\end{document}